\documentclass[]{article}
\usepackage{amsfonts,amsbsy,amscd,amsgen,amsthm,dsfont,amsmath,amssymb,mathrsfs}
\usepackage[colorlinks=true,allcolors=blue]{hyperref}
\usepackage{lipsum}
\usepackage{amsfonts}
\usepackage{graphicx}
\usepackage{epstopdf}
\usepackage{algorithmic}
\usepackage{authblk}
\usepackage[margin=1in]{geometry}
\usepackage{comment}

\usepackage{tikz}
\usetikzlibrary{matrix,shapes,decorations.pathreplacing,fit,backgrounds}

\usepackage{amsopn}

\usepackage{subfigure}
\usepackage{bm} %italic vectors
\theoremstyle{definition}
\usepackage{xcolor,comment}
\usepackage{enumitem}
\usepackage[utf8]{inputenc}
\usepackage[T1]{fontenc}
\newcommand{\id}{\mathrm d}
\newcommand{\vc}{\mathbf}

\renewcommand{\i}{\hat{\mathrm i}}
\renewcommand{\tilde}{\widetilde}
\newcommand{\pard}[2]{\frac{\partial #1}{\partial #2}}

\DeclareMathAlphabet\mathbfcal{OMS}{cmsy}{b}{n}
\DeclareMathOperator*{\argmin}{\text{argmin}}

\newcommand{\R}{\mathbb{R}}
\newcommand{\Rn}{\R^n}

\newcommand{\grad}{\nabla}

\newtheorem{thm}{Theorem}

\newtheorem{rem}{Remark}

\title{Fast and scalable computation of shape-morphing nonlinear solutions with application to evolutional neural networks}
\author{William Anderson}
\author{Mohammad Farazmand\thanks{Corresponding author's email address: \href{mailto:farazmand@ncsu.edu}{farazmand@ncsu.edu}}}
\affil{Department of Mathematics, North Carolina State University,\\ 2311 Stinson Drive, Raleigh, NC 27695-8205, USA}
\date{}

\begin{document}
\maketitle

\begin{abstract}
We develop fast and scalable methods for computing reduced-order nonlinear solutions (RONS).
RONS was recently proposed as a framework for reduced-order modeling of time-dependent partial differential equations (PDEs), where the modes depend nonlinearly on a set of time-varying parameters. RONS uses a set of ordinary differential equations (ODEs) for the parameters to optimally evolve the shape of the modes to adapt to the PDE's solution. This method has already proven extremely effective in tackling challenging problems such as advection-dominated flows and high-dimensional PDEs. However, as the number of parameters grow, integrating the RONS equation and even its formation become computationally prohibitive. Here, we develop three separate methods to address these computational bottlenecks: symbolic RONS, collocation RONS and regularized RONS. We demonstrate the efficacy of these methods on two examples:  Fokker--Planck equation in high dimensions and the Kuramoto--Sivashinsky equation. In both cases, we observe that the proposed methods lead to several orders of magnitude in speedup and accuracy. Our proposed methods extend the applicability of RONS beyond reduced-order modeling by making it possible to use RONS for accurate numerical solution of linear and nonlinear PDEs.
Finally, as a special case of RONS, we discuss its application to problems where the PDE's solution is approximated by a neural network, with the time-dependent parameters being the weights and biases of the network. The RONS equations dictate the optimal evolution of the network's parameters without requiring any training.
\end{abstract}

\section{Introduction}
Anderson and Farazmand \cite{anderson2021} recently proposed reduced-order nonlinear solutions (RONS) as a new framework for deriving reduced-order models for time-dependent PDEs. RONS considers shape-morphing approximate solutions,
\begin{equation}
\hat  u( \vc x, \vc q(t))= \sum_{i=1}^r \alpha_i(t) u_i(\vc x,\pmb\beta_i(t)),
\label{eq:RONS_gen}
\end{equation}
to the PDE which depend nonlinearly on time-varying parameters $\vc q(t) =  \{\alpha_i(t),\pmb\beta_i(t)\}_{i=1}^r$. This is in contrast to most reduced-order models which consider approximate solutions 
$\hat{u} (\vc x, \vc q(t) ) = \sum_{ i  }q_i(t) u_i(\vc x)$ as a linear combination of time-independent modes $u_i$ (see \cite{Willcox2015,rowley2017}, for reviews).
By allowing nonlinear dependence on the parameters, RONS significantly expands the scope of reduced-order modeling, resulting in more accurate reduced models capable of tackling challenging problems such as advection-dominated dynamics.

As we review in Section~\ref{sec:prelims}, RONS uses a set of ordinary differential equations (ODEs) for the optimal evolution of parameters $\vc q(t)\in\mathbb R^n$ by minimizing the instantaneous error between the dynamics of the reduced-order solution $\hat u(\vc x,\vc q(t))$ and the true dynamics of the PDE. Furthermore, RONS ensures that the resulting reduced-order model preserves conserved quantities of the PDE.

There are two main computational bottlenecks that may adversely affect the performance of RONS.
The main computational cost of RONS comes from evaluating the functional inner products which are required to form the reduced-order equations. More specifically, in order to form the RONS reduced-order equations, $\mathcal O(n^2)$ inner products must be evaluated, where $n$ is the number of parameters. Making matters worse, these inner products need to be reevaluated at each time step as the parameter values evolve. 
To compute these inner products, one can use quadrature, Monte Carlo integration, or symbolic computing. As the number of parameters $n$ grows, all these methods become quickly prohibitive. The second computational cost arises from the stiffness of the RONS equations. As mentioned earlier, RONS equations are a set of ODEs for the evolution of parameters $\vc q(t)$. As the number of parameters increases, these ODEs can become stiff and therefore very slow to solve using explicit time integration.

In this paper, we develop three separate methods to address these computational bottlenecks, and thus drastically reduce the computational cost of RONS.
For the first method, which we call symbolic RONS, we assume that the inner products can be computed symbolically. Exploiting the hidden structure of RONS equations, 
we reduce the number of required inner product computations to $\mathcal O(K^2)$ where $K\ll n$ is an integer independent of $n$. Furthermore, because this method uses symbolic computation, the inner products do not need to be recomputed during time stepping. As a result, the computational cost remains low even when the number of parameters $n$ is very large. This scalability allows us to go beyond reduced-order modeling and use RONS as a spectral method where the modes (or basis functions) evolve over time through their nonlinear dependence on time-dependent parameters.

The second method, which we refer to as collocation RONS, introduces a collocation point version of RONS in case symbolic computations are not feasible. 
This method minimizes the discrepancy between the RONS dynamics and the governing PDE on a set of prescribed collocation points.
Collocation RONS is applicable to general nonlinear PDEs, does not require inner product evaluations, and therefore does not require any symbolic computation.
Furthermore, we show that using Monte Carlo integration to approximate the RONS equations coincides with solving a least squares problem which arises from our collocation point method. However, the system of equations arising from collocation RONS is significantly better conditioned than the Monte Carlo approach, and therefore numerically more stable.

Our third contribution addresses the stiffness of the RONS equations. 
When the number of model parameters is large, the RONS ODEs can become stiff and therefore slow to solve using explicit time integration schemes. To address this issue, we introduce a regularized version of RONS which is applicable to both symbolic RONS and collocation RONS.
Regularized RONS introduces a Tikhonov penalization to the underlying minimization problem. This regularization significantly speeds up the numerical time integration of the RONS equations while insignificantly affecting the accuracy of the solutions.

\subsection{Related work}\label{sec:relWork}
Before RONS~\cite{anderson2021}, several previous studies had already considered nonlinear shape-morphing approximate solutions $\hat u(\vc x, \vc q(t)) $ for specific PDEs. For instance, to build reduced-order models for the nonlinear Schr{\"o}dinger (NLS) equation, several authors have proposed approximating wave packets with either Gaussian or hyperbolic secant envelopes~\cite{Adcock12, adcock09,cousins15,PerezGarcia1996,ruban2015,ruban2015b}. 
The amplitude, width, and center of the wave packet are controlled by parameters that evolve over time. Refs.~\cite{PerezGarcia1996,ruban2015,ruban2015b} use the variational Lagrangian formulation of NLS to obtain a set of ODEs for evolving these parameters. As an alternative approach, Adcock et al.~\cite{Adcock12, adcock09} use the symmetries of NLS to evolve the parameters. Another example appears in fluid dynamics where vortex methods approximate the fluid flow as a superposition of point vortices~\cite{newton_Nvortex}, or their smooth approximations~\cite{beale1985,koumoutsakos_2000}. The position, strength, and shape of the vortices are then evolved based on the induced velocity of other vortices.

Although the idea of shape-morphing approximate solutions has been around for decades, the evolution of their shape parameters were determined using ad hoc methods on a case by case basis. RONS proposed a unified framework for evolving these parameters which is applicable to a broad range of PDEs, without relying on the variational structure or symmetries of the PDE.

Interestingly, in the context of evolutional deep neural networks (EDNNs), Du and Zaki \cite{du21} simultaneously and independently derived a set of evolution equations similar to RONS~\cite{anderson2021}. EDNNs approximate solutions of PDEs by evolving weights and biases of a deep neural network over time. Since the network's activation functions are nonlinear, an EDNN depends nonlinearly on its parameters, i.e., weights and biases. As such, EDNNs are a special case of reduced-order nonlinear solutions. Therefore, it is not surprising that the EDNN equations are similar to RONS.

In spite of this similarity, there are some notable differences between EDNN and RONS.
Namely, RONS ensures that the reduced-order model respects the conserved quantities of the PDE.
EDNN does not guarantee these conservation laws, although they can be easily enforced following the methodology introduced in~\cite{anderson2021}.
On the other hand, Du and Zaki \cite{du21} show how various boundary conditions of the PDE can be embedded into the EDNN framework, an important contribution which was not considered in the development of RONS.

As in RONS, forming the EDNN equations requires the evaluation of certain functional inner products.
Du and Zaki~\cite{du21} approximate these inner products using Monte Carlo integration with uniform sampling.
Bruna et al.~\cite{bruna22} proposed an adaptive sampling method to estimate the integrals.
Their adaptive samples are drawn from a distribution which depends on the approximate solution at any given time.
They show that, for PDEs whose solutions are localized in space, adaptive sampling results in more accurate solutions than uniform sampling.

As mentioned earlier, an important feature of RONS is its ability to ensure that the approximate solutions preserve the conserved quantities of the original PDE. There are many studies which consider the same objective; however, the resulting methods are only applicable to a special class of governing equations. For instance, symplectic integrators are 
specifically designed to preserve the two-form associated with a Hamiltonian system~\cite{Bridges2006,McLachlan2006,Cifani2022}.
Similarly, Peng and Mohseni~\cite{peng16} developed proper symplectic decomposition (PSD) for Hamiltonian systems to ensures their reduced-order models preserve the Hamiltonian structure of the full-order model. Carlberg et al.~\cite{carlberg2018} propose a finite-volume based method which guarantees preservation of conserved quantities in the reduced model. This method is only applicable to PDEs derived from conservation laws, and hence amenable to finite-volume discretization. In contrast, RONS preserves any finite number of conserved quantities of the PDE without making any restricting assumptions on the structure of the PDE.

Finally, we point out that the method of optimally time-dependent (OTD) modes~\cite{otd,babaee17,PRE2016,Babaee2022} uses an expansion similar to Eq.~\eqref{eq:RONS_gen}. However, OTD is only applicable to stability analysis of linear or linearized PDEs. In contrast, RONS is applicable for reduced-order modeling and numerical approximation of general nonlinear PDEs.

\subsection{Outline}
This paper is organized as follows. In section \ref{sec:prelims}, we briefly review the derivation of RONS and its relation to Galerkin-type methods.
Section~\ref{sec:methods} contains our main theoretical results where we develop fast and scalable methods for constructing and solving the RONS equations. Section~\ref{sec:num_ex} contains numerical results demonstrating the application of the proposed methods to two different PDEs. We present our concluding remarks in section~\ref{sec:conc}.

\section{Set-up and preliminaries}
\label{sec:prelims}

In this section, we present a succinct review of RONS. We refer to Ref.~\cite{anderson2021}
for a more detailed discussion.
RONS builds reduced-order models for PDEs of the general form
\begin{equation}
\frac{ \partial u }{ \partial t } = F(u), \quad u( \vc x, 0 ) = u_0(\vc x),
\label{eq:general_pde}
\end{equation}
where $u:D \times \R^+ \rightarrow \R^p, (\vc x,t)\mapsto u(\vc x,t)$ is the solution of the PDE, $D \subseteq \R^d$ is the spatial domain,
$F$ is a potentially nonlinear differential operator, and $u_0$ is the initial condition.
We assume the solution $u(\cdot,t)$ belongs to a Hilbert space $H$ with the inner product $\langle \cdot,\cdot \rangle_H$ and the induced norm $\|\cdot \|_H$. To simplify the exposition, we assume $p=1$ hereafter, i.e., $u(\vc x,t)\in\mathbb R$. Generalization to $p>1$ and to complex-valued functions is straightforward~\cite{anderson22}.

We consider \emph{shape-morphing approximate solutions} $\hat u(\vc x,\vc q (t))$ which depend nonlinearly on a set of time-dependent parameters $\vc q(t)\in\mathbb R^n$.
RONS prescribes a set of ODEs to evolve the parameters $\vc q(t)$ such that the instantaneous error between dynamics of the reduced-order solution $\hat u$ and dynamics of the true PDE is minimized. The instantaneous error is defined by
\begin{equation}
\mathcal J( \vc q, \dot{ \vc q }) = \frac{1}{2} \|\hat u_t-F(\hat u)\|_H^2,
\label{eq:costfunctional}
\end{equation}
which measures the difference between the rate of change of the approximate solution $\hat u_t$ and the rate of change $F(\hat u)$ dictated by the PDE. 
Here $\hat{u}_t$ is shorthand for 
\begin{equation}
	\hat{u}_t (\vc x, \vc q(t)) = \sum_{ i = 1 }^n \frac{ \partial \hat{u} }{ \partial q_i }(\vc x, \vc q(t)) \dot{ q_i }.
\end{equation}
If the PDE has no conserved quantity, the reduced-order equations are obtained by minimizing~\eqref{eq:costfunctional}. However, let's consider the more general case where the PDE has 
$m$ conserved quantities, $I_k: H \to \R$ with $k\in\{1,2,\cdots,m\}$.
Since these quantities are conserved, they must satisfy $I_k(u( \cdot ,t)) = I_k(u_0)$ for all $t\geq 0$.
It is desirable for the reduced-order model to also preserves these conserved quantities, since otherwise the reduced model may exhibit unphysical behavior~\cite{peng16,majda2012}.

To obtain an evolution equation for the parameters $\vc q(t)$, we solve the constraint optimization problem
\begin{align}\label{eq:min_cons}
& \min_{\dot{\vc q} \in \R^n} \mathcal J(\vc q,\dot{\vc q}),\nonumber\\
&\mbox{subject to}\quad I_{ k }(\vc q (t))= I_{ k }(\vc q (0)) , \quad k =1,2,...,m, \quad \forall t\geq 0,
\end{align}
%\begin{equation}\label{eq:constJ}
%\min_{\dot{\vc q} \in \R^n} \mathcal J(\vc q,\dot{\vc q}).
%\end{equation}
where $I_{k}( \vc q(t) )$ is shorthand for $I_k(\hat u( \cdot ,\vc q (t) ) )$.
As shown in~\cite{anderson2021}, the solution to this minimization problem is
\begin{equation}
M(\vc q)\dot{\vc q} =  \vc f  (\vc q) - \sum_{ k = 1 }^{ m } \lambda_k \nabla I_k(\vc q),
\label{eq:qdot_const}
\end{equation}
which we refer to as the RONS equation.
Here $M(\vc q) \in \R^{n \times n}$ is the symmetric positive definite \textit{metric tensor} defined by
\begin{equation}
M_{ij} = \left\langle \pard{\hat u}{q_i}, \pard{\hat u}{q_j}\right\rangle_H,\quad i,j\in\{1,2,\cdots,n\}.
\label{eq:metricT}
\end{equation}
The entries of the right-hand side vector field $\vc f :\Rn\to \Rn$ are given by
\begin{equation}
f_i = \left\langle \pard{\hat u}{q_i}, F(\hat u)\right\rangle_H,\quad i=1,2,\cdots, n.
\label{eq:rhs}
\end{equation}
The Lagrange multipliers $\pmb \lambda =  ( \lambda_1, ... , \lambda_m )^\top$ satisfy the linear system
\begin{equation}
C(\vc q)  \pmb \lambda = \vc b(\vc q) ,
\label{eq:lambda_eq}
\end{equation}
where $C(\vc q)\in \R^{m\times m}$ is the symmetric positive definite \textit{constraint matrix} defined by
\begin{equation}
C_{ij} = \langle \grad I_j ,M^{-1} \grad I_i\rangle,\quad i,j\in\{1,2,\cdots,m\}, 
\label{eq:constMat}
\end{equation}
and the vector $\vc b= (b_1,b_2,\cdots,b_m)^\top\in\R^m$ is given by
\begin{equation}
b_i = \langle \grad I_i ,M^{-1} \vc f \rangle,\quad i=1,2,\cdots, m,
\label{eq:b}
\end{equation}
where $\langle \cdot, \cdot \rangle$ denotes the standard Euclidean inner product. The gradients $\nabla I_i$ denote the partial derivatives with respect to the components of the parameters $\vc q$.

If no conserved quantities are enforced, then we must solve the optimization problem~\eqref{eq:min_cons} without any constraints. 
The unique minimizer of the unconstrained problem is given by omitting the summation term from Eq.~\eqref{eq:qdot_const}, i.e.,
\begin{equation}
M(\vc q)\dot{\vc q} = \vc f  (\vc q).
\label{eq:qdot_unconst}
\end{equation}
As we mentioned in section~\ref{sec:relWork}, equation~\eqref{eq:qdot_unconst} was derived simultaneously and independently by Du and Zaki~\cite{du21} in the context of EDNNs. However, the more general equation~\eqref{eq:qdot_const}, which ensures the preservation of conserved quantities, was only derived in Ref.~\cite{anderson2021}.

\begin{figure*}
	\centering
	\includegraphics[width=0.95\textwidth]{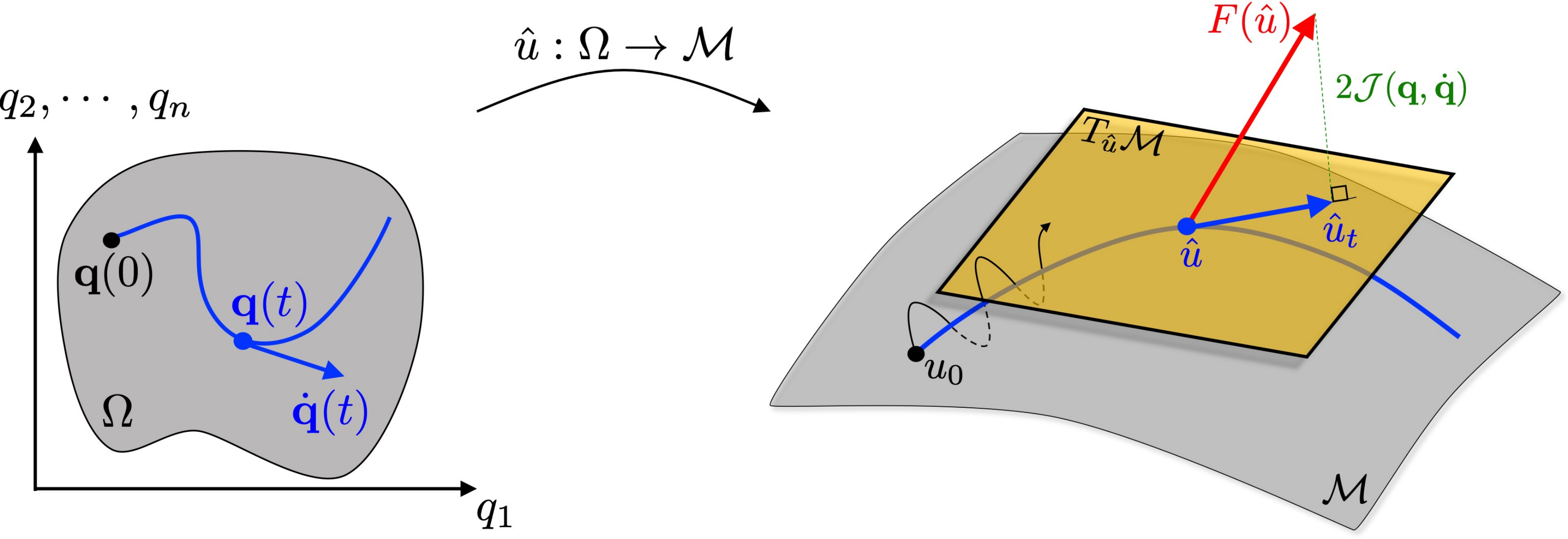}
	\caption{Geometric illustration of RONS. The shape-morphing approximate solution $\hat u$ maps the parameter space $\Omega$ to the manifold $\mathcal M$. An evolution of the parameters $\vc q(t)$ is mapped to a curve on $\mathcal M$, and the tangent vector to the parameters $\dot{\vc q}(t)$ is mapped to $\hat u_t$ in the tangent space of the manifold. }
	\label{fig:Geometric_Explanation}
\end{figure*}

A geometric depiction of RONS in the unconstrained case is shown in Figure \ref{fig:Geometric_Explanation}. 
We view the shape-morphing approximate solution $\hat u(\cdot ,\vc q)$ as a map from the parameters $\vc q$ to the Hilbert space $H$ where solutions of the PDE lie.
The approximate solution $\hat u$ maps the set of all viable parameter values $\vc q\in\Omega \subseteq \R^n$ to an $n$-dimensional manifold $\mathcal{M} \subset H$.
An arbitrary but smooth evolution of parameters $\vc q(t)$ defines a smooth curve in the set $\Omega$. The tangent vector, or velocity, of this curve is given by $\dot{\vc q} (t)$. Under the map $\hat u$, this curve
is mapped onto a curve which lies on the manifold $\mathcal M$ in the function space $H$. The tangent vector $\dot{\vc q} (t)$ is mapped to the tangent vector $\hat u_t$ which lies on the tangent space of the manifold $T_{\hat u} \mathcal M$.
In general, the manifold is not invariant under the dynamics of the governing PDE~\eqref{eq:general_pde}, and therefore $F(\hat u)$ will not necessarily lie in tangent space. By minimizing~\eqref{eq:costfunctional} with respect to $\dot{\vc q}(t)$, we find the vector $\hat u_t$ which is the orthogonal projection of $F(\hat u)$ onto $T_{\hat u} \mathcal M$. 
In other words, we evolve the approximate solution $\hat u$ so that it most closely resembles the expected PDE dynamics. 

Conventional Galerkin projection models are a special case of RONS. Consider an approximate solution which depends linearly on the parameters,
\begin{equation}
\hat u(\vc x, \vc q (t)) = \sum_{i=1}^{n} q_i(t) u_i (\vc x), 
\label{eq:uhat_galerkin}
\end{equation}
where the modes $\{ u_i\}_{i=1}^n$ are prescribed orthonormal functions, e.g., proper orthogonal decomposition (POD) modes. In this special case, the unconstrained RONS equation~\eqref{eq:qdot_unconst} coincides with standard Galerkin projection.
More specifically, the metric tensor $M(\vc q)$ becomes the identity matrix and the right-hand side vector is given by $f_i = \langle u_i,F(\hat u)\rangle_H$. From a geometric standpoint, for the reduced-order solution~\eqref{eq:uhat_galerkin}, the manifold $\mathcal M$ becomes a flat subspace spanned by the modes $\{ u_i\}_{i=1}^n$. We refer to~\cite{anderson2021} for further details.

\section{Fast and scalable computational methods}
\label{sec:methods}
\subsection{Computational bottlenecks}
Although RONS has shown great promise for both reduced-order modeling and numerical simulation of PDEs~\cite{anderson2021,du21,bruna22,anderson22}, forming and solving the RONS equations can be computationally expensive. In this section, we first outline the main computational bottlenecks associated with RONS and then present our proposed remedies.

The computational cost of RONS equations~\eqref{eq:qdot_const} comes from three main sources:
\begin{enumerate}
	\item Forming the metric tensor $M(\vc q)$ and the vector field $\vc f(\vc q)$ (addressed in sections~\ref{sec:sym_comp} and~\ref{sec:colloc}).
	\item Stiffness of the RONS equations (addressed in section~\ref{sec:tikh}).
	\item Inverting the metric tensor. By inverting the metric tensor we refer to any numerical method for solving the linear equation~\eqref{eq:qdot_const} for $\dot{\vc q}$.
\end{enumerate}

We now describe each of these computational bottlenecks in more detail. To form the metric tensor~\eqref{eq:metricT}, we need to compute $n^2$ inner products. Since this matrix is symmetric, the number of independent inner products is in fact $n(n+1)/2$. Additionally, to form the right-hand side vector~\eqref{eq:rhs}, we need to compute $n$ inner products. Therefore, a total of $n(n+3)/2$ integrals need to be computed. As the number of parameters $n$ grows, this becomes computationally prohibitive. Making matters worse, during time stepping, $\vc q(t)$ changes and these integrals need to be recomputed at each time step. In section~\ref{sec:sym_comp}, we develop a method which drastically reduces the number of inner product computations. We refer to this method as \emph{symbolic RONS}, or \emph{S-RONS} for short.

Symbolic RONS requires the inner products to be symbolically computable. Depending on the choice of the approximate solution $\hat u$, this may not be feasible. In section~\ref{sec:colloc}, we develop a collocation point approach to  RONS which does not require any integral or inner product computation and therefore reduces the computational cost of RONS by several orders of magnitude. We refer to this method as \emph{collocation RONS}, or \emph{C-RONS} for short.

The second issue arises from the fact that the RONS equations~\eqref{eq:qdot_const} can be stiff as a set of ODEs. As a result, using explicit schemes for time integration may require exceedingly small time steps. In section~\ref{sec:tikh}, we propose a regularized version of the optimization problem~\eqref{eq:min_cons} which alleviate this issue. We refer to the resulting method as the \emph{regularized RONS}, which can be used in conjunction with both S-RONS and C-RONS.

In our experience, the last issue (inverting the metric tensor) does not present a major roadblock. There exist several fast methods for solving large linear systems which can be used for inverting the metric tensor~\cite{Gill2021}. Therefore, we focus on items 1 and 2 above which constitutes the main computational bottlenecks.

\subsection{Symbolic RONS}\label{sec:sym_comp}
In this section, we present a method for efficient construction of the metric tensor $M$ and right-hand side vector $\vc f$ using symbolic computation of the required inner products.
We show that only a relatively small number of symbolic computations are required to build the RONS equation, provided that the shape-morphing approximate solution has a specific form and analytical symbolic expressions for the inner products in $M$ and $\vc f$ can be obtained.

Specifically, we consider shape-morphing approximations of the form,
\begin{equation}
\hat u (\vc  x, \vc q (t) ) = \sum_{ i = 1 }^{r} \alpha_i (t) \phi(\vc x, \pmb \beta_i (t) ),
\label{eq:ansatz}
\end{equation}
where $\phi(\cdot,\pmb\beta_i):\mathbb R^d\to \mathbb R$ is a $C^1$ function. We refer to $\phi(\vc x, \pmb \beta_i (t) )$ as the $i$-th \textit{shape-morphing mode}. The vector of \textit{shape parameters} $\pmb \beta_i \in \R^{K-1}$ controls the shape of the $i$-th mode; it contains parameters such as length scales and center of the mode. The scalar $\alpha_i (t)$ denotes the mode amplitude. Therefore, parameters of the shape-morphing solution $\hat u$ are given by $\vc q = (\alpha_1, \pmb \beta_1^\top, ..., \alpha_r, \pmb \beta_r^\top)^\top \in  \R^{rK}$, where $n=rK$. Note that the number of shape parameters for each mode $K-1$ is independent of the number of terms $r$ in the sum. This independence plays an important role in the proposed computational method.

As an example, we can consider Gaussian modes $\phi$ which lead to the approximate solution,
\begin{equation}
\hat u( x, \vc  q (t)) = \sum_{ i = 1 }^{r} A_i(t) \exp \bigg[- \frac{ (  x -   c_i(t) )^2 }{ L_i(t)^2  } \bigg],
\label{eq:gauss_ans}
\end{equation}
where $\alpha_i(t) = A_i(t)$ and $\pmb \beta_i(t) = (L_i(t),  c_i(t))^\top$. 
Here $A_i$ controls the amplitude of the $i$th Gaussian, $L_i$ is a length scale that determine the Gaussian's width, and $ c_i$ determines the Gaussian's center. We use this one-dimensional Gaussian mixture as an illustrative example throughout this section.

There are many other possible choices of modes $\phi$ for the approximate solution. 
For example, we could take the modes to be activation functions typically used in neural networks, such as the rectified linear unit (ReLU) or hyperbolic tangent. In this case, Eq.~\eqref{eq:ansatz} represents a shallow neural network and the shape parameters $\pmb\beta_i$ are the weights and biases of the $i$-th node. Another choice could be wavelet functions, where the shape parameters are dilations and translations.

We now exploit the structure of the shape-morphing approximation~\eqref{eq:ansatz} to efficiently calculate the inner products in $M$ and $\vc f$ using symbolic computation. We first discuss how to efficiently build the metric tensor $M$.
Consider arbitrary indices $i,j \in \{1, ... ,r\}$.
Using these general indices and the approximate solution~\eqref{eq:ansatz}, all entries of $M$ will have the form of one of the following inner products,
\begin{subequations}\label{eq:ips_M}
	\begin{align}
	\mathcal{I}_{\alpha_i \alpha_j} &:=\bigg\langle \frac{ \partial \hat u}{ \partial \alpha_i } , \frac{ \partial \hat u}{ \partial \alpha_j } \bigg\rangle_{H},&  \\
	\mathcal{I}_{\alpha_i \beta_{jk}} &:= \bigg\langle \frac{ \partial \hat u}{ \partial \alpha_i } , \frac{ \partial \hat u}{ \partial \beta_{jk} } \bigg\rangle_{H},& \quad &k \in \{ 1, ..., K-1 \}& \\
	\mathcal{I}_{\beta_{ik} \beta_{j\ell}} &:= \bigg\langle \frac{ \partial \hat u}{ \partial \beta_{ik} } , \frac{ \partial \hat u}{ \partial \beta_{j\ell} } \bigg\rangle_{H},& \quad &k, \ell \in \{ 1, ..., K-1 \}&,
	\label{eq:beta_derivs}
	\end{align}
\end{subequations}
where $\beta_{ik}$ denotes the $k$-th component of $\pmb \beta_i=(\beta_{i1},\beta_{i2},\cdots, \beta_{i(K-1)})^\top$.
The advantage of using symbolic computation for the expressions in equation~\eqref{eq:ips_M} is that after obtaining closed-form expressions for the inner products, we can build the entire metric tensor through substitution of the appropriate indices $i$ and $j$.
For example, rather than computing $\mathcal{I}_{\alpha_1 \alpha_2}$, we simply need to substitute the values of $\alpha_1$ and  $\alpha_2$ into the already obtained symbolic expression for $\mathcal{I}_{\alpha_i \alpha_j}$. This same principle holds regardless of which indices we choose as $i$ and $j$.

Note that, after obtaining closed-form expressions the inner products in equation~\eqref{eq:ips_M}, we can evaluate all entries of the metric tensor $M$ regardless of the number of modes $r$ used in the approximate solution.
Additionally, we only need to perform the symbolic computations at the initial time and can then substitute the updated parameter values as we march the approximate solution foward in time.
Note that, by symmetry of the inner product, we only need to calculate $K(K-1)/2$ inner products for equation~\eqref{eq:beta_derivs}. Therefore, there are in total $K(K+1)/2$ terms to be calculated in equation~\eqref{eq:ips_M}. We emphasize that this number is independent of the number of modes $r$ in the shape-morphing approximation~\eqref{eq:ansatz}; it only depends on the number of shape parameters $K$.

Symbolic RONS (or S-RONS) can be alternatively described by examining the structure of the metric tensor.  The matrix $M$ is composed of blocks $M^{(i,j)} \in \R^{ K\times K} $ such that
\begin{equation}
M = \begin{bmatrix}
M^{(1,1)} & M^{(1,2)} & \cdots  & M^{(1,r)}  \\
M^{(2,1)} & M^{(2,2)} &\cdots & M^{(2,r)} \\
\vdots & \vdots & \ddots &   \vdots\\
M^{(r,1)} & M^{(r,2)} &  \cdots & M^{(r,r)}
\end{bmatrix},
\end{equation}
where $M^{(i,j)}= (M^{(j,i)})^\top$ since $M$ is symmetric.
Each block can be expressed in terms of the inner products~\eqref{eq:ips_M}, 

\begin{equation} \label{eq:M_block}
 \begin{tikzpicture}[baseline,
%Styles
Matrix/.style={
	matrix of nodes,
	text height=2.5ex,
	text depth=0.75ex,
	text width=15ex,
	align=center,
	left delimiter={[},
	right delimiter={]},
	column sep=0pt,
	row sep=0pt,
	%nodes={draw=black!10}, % Uncoment to see the square nodes.
	nodes in empty cells,
},
DA/.style={
	fill=green!10,
	solid,
	line width=0.5pt,
	rounded corners=0pt,
}
]
\node at (-7.5,0) {$M^{(i,j)} =$};
\matrix[Matrix] at (0,0) (M2){  
	$\mathcal{I}_{\alpha_i \alpha_j} $ &  $\mathcal{I}_{ \alpha_i \beta_{j1}  } $ &   $ \cdots   $ &  $\mathcal{I}_{ \alpha_i \beta_{j(K-1)} }  \medskip $\\ 
	$\mathcal{I}_{\beta_{i1} \alpha_j  } $ & $\mathcal{I}_{\beta_{i1} \beta_{j1}} $&  $ \cdots  $ & $  \mathcal{I}_{ \beta_{i1} \beta_{j(K-1)}  } \medskip $\\ 
	$\vdots $&  $\vdots $  & $\ddots$ &   $\vdots \medskip $\\
	$\mathcal{I}_{ \beta_{i(K-1)} \alpha_j } $ &$\mathcal{I}_{ \beta_{i(K-1)}\beta_{j1}  }$ &  $\cdots$ &$ \mathcal{I}_{\beta_{i(K-1)}\beta_{j(K-1)}} $\\
};
\begin{scope}[on background layer] %Allows to draw in background layer.
%FOR MATRIX M
%To delimit area groups using matrix nodes
% Correcting using the method you use 
\node[fit=(M2-1-2)](n1){};
\node[fit=(M2-2-3)](n2){};
\node[fit=(M2-3-4)](n3){};  
%\node[fit=(M2-4-4)](n4){};
%If you trace the draw:
\draw[DA,fill=none,fill opacity=0.1](n1.north west)
-| (n3.south east)
%|- (n4.south west)
|- (n3.south west)
|- (n2.south west)
|- (n1.south west)
-- cycle;           
\end{scope}
\end{tikzpicture}
.
\end{equation}

%\begin{equation}
%M^{(i,j)} = \begin{bmatrix}
%\mathcal{I}_{\alpha_i \alpha_j} & \mathcal{I}_{ \alpha_i \beta_{j1}  }  & \cdots & \cdots  & \mathcal{I}_{ \alpha_i \beta_{j(K-1)}  } \medskip \\ 
%\mathcal{I}_{\beta_{i1} \alpha_j  }  & \mathcal{I}_{\beta_{i1} \beta_{j1}} & \cdots & \cdots & \mathcal{I}_{ \beta_{i1} \beta_{j(K-1)}  } \medskip \\ 
%\vdots &  \vdots & \mathcal{I}_{\beta_{i2}\beta_{j2}} &  & \vdots \\ 
%\vdots &  \vdots &  \vdots  & \ddots & \vdots\medskip \\
%\mathcal{I}_{ \beta_{i(K-1)} \alpha_j }  &\mathcal{I}_{ \beta_{i(K-1)}\beta_{j1}  } & \mathcal{I}_{ \beta_{i(K-1)} \beta_{j2} }  & \cdots & \mathcal{I}_{\beta_{i(K-1)}\beta_{j(K-1)}} 
%\end{bmatrix}.
%\end{equation}  
With this labeling, the block $M^{(i,j)}$ represents all of the inner products involving derivatives of the approximate solution with respect to the shape parameters $\beta_{ik}$ and $\beta_{jk}$ of the $i$-th and $j$-th modes and their respective amplitudes, $\alpha_i$ and $\alpha_j$.
Although there are $K^2$ entries in $M^{(i,j)}$, we only need to perform symbolic calculations for the $K(K+1)/2$ entries in the lower triangular part of the block. The remaining entries of the matrix, enclosed in a box in~\eqref{eq:M_block}, are then determined by the symmetry of inner products in~\eqref{eq:ips_M}.
In other words, if we have a closed-form expression for an entry $\mathcal{I}_{q_i q_{j}  }$ in the lower triangular part of $M^{(i,j)}$, we obtain the corresponding entry $\mathcal{I}_{q_j q_{i } }$ in the upper triangular part of the block by simply swapping the values of $q_i$ and $q_j$ in the symbolic expression.

As an example, consider the Gaussian mixture~\eqref{eq:gauss_ans}. 
We must symbolically compute six inner products to form the block,

\begin{equation}\label{eq:M_ij_gm}
\begin{tikzpicture}[baseline,
%Styles
Matrix/.style={
	matrix of nodes,
	text height=4ex,
	text depth=2ex,
	text width=15ex,
	align=center,
	left delimiter={[},
	right delimiter={]},
	column sep=0pt,
	row sep=0pt,
	%nodes={draw=black!10}, % Uncoment to see the square nodes.
	nodes in empty cells,
},
DA/.style={
	fill=green!10,
	solid,
	line width=0.5pt,
	rounded corners=0pt,
}
]
\node at (-6,0) {$M^{(i,j)} =$};
\matrix[Matrix] at (0,0) (M2){ % Matrix contents  
	$\bigg\langle \dfrac{ \partial \hat u}{ \partial A_i } , \dfrac{ \partial \hat u}{ \partial A_j } \bigg\rangle_{H} $& $ \bigg\langle \dfrac{ \partial \hat u}{ \partial A_i } , \dfrac{ \partial \hat u}{ \partial L_j } \bigg\rangle_{H} $& $ \bigg\langle \dfrac{ \partial \hat u}{ \partial A_i } , \dfrac{ \partial \hat u}{ \partial c_j } \bigg\rangle_{H} \medskip $\\ 
	$\bigg\langle \dfrac{ \partial \hat u}{ \partial L_i } , \dfrac{ \partial \hat u}{ \partial A_j } \bigg\rangle_{H} $& $\bigg\langle \dfrac{ \partial \hat u}{ \partial L_i } , \dfrac{ \partial \hat u}{ \partial L_j } \bigg\rangle_{H} $& $ \bigg\langle \dfrac{ \partial \hat u}{ \partial L_i } , \dfrac{ \partial \hat u}{ \partial c_j } \bigg\rangle_{H} \medskip$ \\ 
	$\bigg\langle \dfrac{ \partial \hat u}{ \partial c_i } , \dfrac{ \partial \hat u}{ \partial A_j } \bigg\rangle_{H} $&  $\bigg\langle \dfrac{ \partial \hat u}{ \partial c_i } , \dfrac{ \partial \hat u}{ \partial L_j } \bigg\rangle_{H} $& $\bigg\langle \dfrac{ \partial \hat u}{ \partial c_i } , \dfrac{ \partial \hat u}{ \partial c_j } \bigg\rangle_{H}$\\
};
\begin{scope}[on background layer] %Allows to draw in background layer.
%FOR MATRIX M
%To delimit area groups using matrix nodes

% Correcting using the method you use 
\node[fit=(M2-1-2)](n1){};
\node[fit=(M2-2-3)](n2){};
\node[fit=(M2-3-4)](n3){};  
%\node[fit=(M2-4-4)](n4){};
%If you trace the draw:
\draw[DA,fill=none,fill opacity=0.1](n1.north west)
-| (n2.south east)
%|- (n4.south west)
|- (n2.south west)
%|- (n2.south west)
|- (n1.south west)
-- cycle;           
\end{scope}
\end{tikzpicture}.
\end{equation}
Note that the terms enclosed in the box can be evaluated using the lower triangular part of the matrix. For instance, $\langle \partial_{A_i}\hat u,\partial_{L_j}\hat u\rangle$
is evaluated using the symbolic expression for $\langle \partial_{L_i}\hat u,\partial_{A_j}\hat u\rangle$ by substituting the values of $A_i$ and $L_j$ instead of $A_j$ and $L_i$, respectively. Therefore, only 6 symbolic computations are required to form the matrix block~\eqref{eq:M_ij_gm} and consequently the entire metric tensor $M$. In comparison, computing the metric tensor by a brute force method, such as quadrature or Monte Carlo methods, would require evaluating $n(n+1)/2=3r(3r+1)/2$ integrals, which becomes prohibitive as the number of terms $r$ increases.

The idea for building $\vc f$ is similar to that of the metric tensor.
We again consider a general index $i\in \{1,...,r\}$ and note all entries of $\vc f$ will have the form of one of the following inner products:
\begin{equation}
	\bigg\langle \frac{ \partial \hat u}{ \partial \alpha_i } , F(\hat u) \bigg\rangle_{H}, \quad
	\bigg\langle \frac{ \partial \hat u}{ \partial \beta_{ik} } , F(\hat u) \bigg\rangle_{H}, \quad k = 1, ..., K-1.
	\label{eq:ips_f}
\end{equation}
After using symbolic computation to obtain closed-form expressions for the $K$ inner products in equation~\eqref{eq:ips_f}, we can then build $\vc f$ through substitution rather than individually calculating each of the $n=rK$ inner products in $\vc f$.

The vector field $\vc f$ also has a block structure.
We may consider $\vc f$ as $r$ vectors $ \vc f^{ (i) } \in \R^{ K }$ stacked on top of each other so that
\begin{equation}
\vc f = \begin{bmatrix}
\vc f^{ (1) }\\ 
\vdots \\
\vc f^{ (r) }
\end{bmatrix},
\end{equation}
where each vector $\vc f^{ (i) }$ is defined by
\begin{equation}
\vc f^{ (i) } = \begin{bmatrix}
\left\langle \dfrac{ \partial \hat u}{ \partial \alpha_i}, F(\hat u)\right\rangle_H, & 
\left\langle \dfrac{ \partial \hat u}{ \partial \beta_{i1}}, F(\hat u)\right\rangle_H , & 
\cdots, &
\left\langle \dfrac{ \partial \hat u}{ \partial \beta_{i(K-1)} }, F(\hat u)\right\rangle_H 
\end{bmatrix}^\top.
\end{equation}
%\begin{equation}
%	\vc f^{ (i) } = \begin{bmatrix}
%	\left\langle \dfrac{ \partial \hat u}{ \partial \alpha_i}, F(\hat u)\right\rangle_H \medskip\\ 
%	\left\langle \dfrac{ \partial \hat u}{ \partial \beta_{i1}}, F(\hat u)\right\rangle_H \medskip\\
%	\vdots \\
%	\left\langle \dfrac{ \partial \hat u}{ \partial \beta_{i2} }, F(\hat u)\right\rangle_H 
%	\end{bmatrix}.
%\end{equation}
To evaluate the entire vector $\vc f$, we only need the symbolic expression for one of the blocks $\vc f^{(i)}$.
Again using the Gaussian mixture as an example, we have
\begin{equation}
\vc f^{ (i) } = \begin{bmatrix}
\left\langle \dfrac{ \partial \hat u}{ \partial A_i}, F(\hat u)\right\rangle_H , & 
\left\langle \dfrac{ \partial \hat u}{ \partial L_i}, F(\hat u)\right\rangle_H , & 
\left\langle \dfrac{ \partial \hat u}{ \partial c_i }, F(\hat u)\right\rangle_H 
\end{bmatrix} ^\top .
\end{equation}
Using a general index $i$ we need symbolic expressions for only three integrals to build the vector $\vc f$ through substitution rather than computing $3r$ integrals. 

In summary, building the metric tensor $M$ requires $K(K+1)/2$ symbolic integrations and building the right-hand side vector $\vc  f$ requires $K$ symbolic computations, resulting in only $K(K+3)/2$ symbolic computations to evaluate $n(n+1)$ terms appearing in the RONS equation~\eqref{eq:qdot_unconst}.
The above discussion leads to the following theorem.
\begin{thm} \label{thm:comptrick}
	Consider a shape-morphing approximate solution of the form~\eqref{eq:ansatz}. Forming the metric tensor $M$ and the right-hand side vector $\vc f$ in the RONS equation~\eqref{eq:qdot_const} requires symbolic calculation of $K(K+3)/2$ inner products. The number of symbolic computations is independent of the number of modes $r$ used in the approximate solution. 
\end{thm}

An important implication of Theorem~\ref{thm:comptrick} is that the number of modes $r$ can be increased arbitrarily without making the computations prohibitive. As a result, it extends RONS beyond a reduced-order modeling framework, where relatively small number of modes are used, and allows us to use RONS for accurate approximation of the PDE's solutions. More specifically, one can think of~\eqref{eq:ansatz} as a spectral method with the modes $\phi(\cdot,\pmb\beta_i(t))$. In contrast to conventional spectral methods, such as the Fourier spectral method, the modes are allowed to change their shape and position over time to adapt to the solution of the PDE by evolving the shape parameters $\pmb\beta_i(t)$. As we show in section~\ref{sec:num_ex}, this shape-morphing property is specially appealing for advection-dominated problems or high-dimensional PDEs with localized solutions.

Finally, we point out that the summation term in equation~\eqref{eq:qdot_const} involves the Lagrange multipliers $\lambda_k$ which are obtained as the solution to the linear system~\eqref{eq:lambda_eq}. Note that this linear system only involves Euclidean inner products and therefore its construction is not computationally expensive.

\begin{comment}
With this theorem,  we have a tool to quickly form the RONS equation with an arbitrary number of modes $r$.
This is particularly useful when modeling systems with complicated dynamics which require many modes to properly capture behavior of true solutions to the governing PDE.
\end{comment}

\subsection{Collocation RONS}
\label{sec:colloc}

While S-RONS is efficient, obtaining symbolic expressions for the required inner products may not always be feasible.
In this section, we present a new approach to RONS where we only enforce that the approximate solution satisfies the governing PDE on a set of prescribed collocation points. This method is applicable for any choice of the approximate solution $\hat u$ and does not require symbolic computing or numerical integration.

We first describe collocation RONS without enforcing any conserved quantities. We define the residual function,
\begin{equation}
	R(\vc x, \vc q, \dot{ \vc q} ) := \hat{u}_t - F(\hat u) = \sum_{ j = 1 }^n \frac{\partial \hat u}{\partial q_j} \dot{q}_j - F(\hat u).
	\label{eq:residual}
\end{equation}
The residual function $R$ measures the point-wise difference between the rate of change of the approximate solution $\hat u_t$ and the dynamics $F(\hat u)$ dictated by the governing PDE. Previous studies \cite{anderson2021,anderson22,bruna22,du21} have all sought an evolution of parameters $\vc q(t)$ which minimizes the norm of the residual function in the underlying Hilbert space by minimizing the cost function~\eqref{eq:costfunctional}.

Here we propose a different approach. In its most general form, we assume $R(\cdot, \vc q,\dot{\vc q}):D\to\mathbb R$ belongs to a function space $V$ to be specified shortly. For any test function $\phi\in V^\ast$, we require $\langle  \phi  , R  \rangle =0$, where $\langle \cdot,\cdot\rangle$ denotes the natural pairing between $V$ and its dual space $V^\ast$.
For computational purposes, we reduce this problem to its finite-dimensional version. More specifically, as in the Petrov-Galerkin method~\cite{karniadakis2005}, we choose a finite number of test functions $\{\phi_i\}_{i = 1}^{N}\in V^\ast$ and require that $\langle  \phi_i , R  \rangle=0$ for $i=1,2,\cdots, N$.

As a special case, we derive a collocation method by assuming that $R(\cdot,\vc q,\dot{\vc q})$ is bounded, i.e., $V=L^\infty(D)$. Furthermore, we consider the test functions $\phi_i (\vc x)= \delta (\vc x - \vc x_i)$ for $ i = 1,2, ...,  N$, where each $\vc x_i$ is a collocation point in the spatial domain $D$.
Note that $\phi_i \in L^1(D)$ and the natural pairing implies 
\begin{equation}
\langle \phi_i,R\rangle = \int_D \delta (\vc x-\vc x_i) R(\vc x,\vc q,\dot{\vc q})\id \vc x = R(\vc x_i,\vc q,\dot{\vc q})=0,
\label{eq:res_coll}
\end{equation}
which requires the residual function to vanish at the collocation point $\vc x_i$.
Using definition~\eqref{eq:residual}, we obtain
\begin{equation}
	\sum_{ j = 1 }^N \frac{\partial \hat u}{\partial q_j}(\vc x_i,\vc q) \dot{q}_j = F(\hat u) \bigg|_{\vc x = \vc x_i}, \quad  i \in \{1,2,\cdots,N\}.
	\label{eq:weak_dirac}
\end{equation}
We write~\eqref{eq:weak_dirac} as a system of equations,
\begin{equation}
\tilde{M}(\vc q) \dot{ \vc q } = \tilde{ \vc f} (\vc q),
\label{eq:rons_colloc}
\end{equation}
where the collocation matrix $\tilde{M}(\vc q)\in \R^{N\times n}$ is given by
\begin{equation}
	\tilde{M}_{ij}(\vc q) = \frac{\partial \hat u}{\partial q_j} (\vc x_i, \vc q), \quad i \in \{1,2,\cdots,N\},  \quad j  \in \{1,2,\cdots,n\} ,
	\label{eq:M_coll}
\end{equation}
and the vector $\tilde{ \vc f}(\vc q)\in \R^{N}$ is defined by
\begin{equation}
\tilde{f}_{i} (\vc q) =  F(\hat u(\vc x, \vc q)) \bigg|_{\vc x = \vc x_i}, \quad i \in \{1,2,\cdots,N\}.
\label{eq:f_coll}
\end{equation}
Note that equation~\eqref{eq:rons_colloc} requires that the approximate solution $\hat u$ satisfies the governing PDE at the collocation points $\vc x_i$.

We refer to equation~\eqref{eq:rons_colloc} as collection RONS, or C-RONS for short. Although this equation resembles the unconstrained RONS equation~\eqref{eq:qdot_unconst}, there are notable differences.
First, to form the C-RONS equation, numerical integrations or symbolic computations are not required; we only need point-wise evaluation of known functions in~\eqref{eq:M_coll} and~\eqref{eq:f_coll}. Second, unlike the metric tensor $M(\vc q)\in\mathbb R^{n\times n}$, the collocation matrix $\tilde M(\vc q)$ is rectangular. Consequently, the linear system~\eqref{eq:rons_colloc} may not have a unique solution.

If the number of collocation points is greater than the number of parameters, $N>n$, then the system is overdetermined and a solution may not exist. 
If the number of collocation points is less than the number of parameters, $N<n$, then the system is underdetermined and there may exist infinitely many solutions to the problem. 
In either case, we obtain $\dot{\vc q}$ by solving the least squares problem,
\begin{equation}
    \min_{\dot{\vc q} \in \R^n} \|  \tilde{M} \dot{\vc q}  - \tilde{ \vc f} \|_2^2,
	\label{eq:coll_problem}
\end{equation}
using the Moore-Penrose pseudoinverse of $\tilde{M}$.
Thus, for collocation RONS, the evolution of parameters is given by
\begin{equation}
	\dot{ \vc q } = \tilde{M}^+(\vc q) \tilde{ \vc f} (\vc q),
	\label{eq:coll_solution}
\end{equation}
where $\tilde M^+(\vc q)\in\mathbb R^{n\times N}$ denotes the pseudoinverse of the collocation matrix $\tilde M(\vc q)$.
If the C-RONS equation~\eqref{eq:rons_colloc} is overdetermined, then the solution~\eqref{eq:coll_solution} is the unique solution to the least-squares problem~\eqref{eq:coll_problem}.
If the system of equations is underdetermined, then the solution~\eqref{eq:coll_solution} is the solution to the least-squares problem with minimal Euclidean norm~\cite{Ipsen2009}.

Finally, we note that the least square problem~\eqref{eq:coll_problem} is equivalent to minimizing the residual sum $\sum_{i=1}^N |R(\vc x_i,\vc q,\dot{\vc q})|^2$ over all possible $\dot{\vc q}\in\mathbb R^n$. In other words, instead of requiring the residual function to vanish at the collocation points as in~\eqref{eq:res_coll}, we choose $\dot{\vc q}$ so that the sum of squares of the residual is minimized.

Now we turn to the problem of enforcing the PDE's conserved quantities in the approximate solution. Note from section~\ref{sec:prelims}, that the governing PDE may have a number of conserved quantities $I_i$ for $i=1,2,\cdots, m$. We would like to ensure that these quantities
are also conserved along the approximate solution $\hat u(\cdot,\vc q(t))$. In other words, we require $I_i(\vc q(t)) = I_{i}(\vc q(0))$ for all times $t\geq 0$.
Taking the derivative of this identity with respect to time, we obtain the equivalent set of equations,
\begin{equation}
	\langle \nabla I_i(\vc q),\dot{\vc q}\rangle = 0, \quad i=1,2,\cdots, m.
\end{equation}

These constraints, together with the C-RONS equation~\eqref{eq:rons_colloc}, lead to the larger system of equations, 
\begin{equation}
	\tilde M_c(\vc q)\dot{\vc q} = \vc{\tilde f}_c(\vc q), 
	\label{eq:rons_colloc_const}
\end{equation}
where the constrained collocation matrix $\tilde M_c\in\mathbb R^{(N+m)\times n}$ and the constrained vector field $\vc{\tilde f}_c\in \mathbb R^{N+m}$ are defined by
\begin{equation}
	\tilde M_c(\vc q) := \begin{bmatrix}
		\tilde M(\vc q) \\
		\nabla I_1(\vc q)^\top \\
		\nabla I_2(\vc q)^\top\\
		\vdots\\
		\nabla I_m(\vc q)^\top
	\end{bmatrix},\quad 
	\vc{\tilde f}_c(\vc q) := \begin{bmatrix}
		\vc{\tilde f}(\vc q) \\
		0 \\
		0\\
		\vdots\\
		0
	\end{bmatrix}.
\end{equation}
As before, we solve the linear system~\eqref{eq:rons_colloc_const} using the pseudoinverse to obtain the constrained collocation equation
$\dot{\vc q} = \tilde M_c^+ \tilde{\vc f}_c$.

\subsubsection{Relation between Monte Carlo sampling and C-RONS}
Monte Carlo integration has previously been used to approximate the inner products in the RONS equation \cite{bruna22,du21}.
We will show that, under certain conditions, the evolution of parameters provided by Monte Carlo approximation coincides with C-RONS~\eqref{eq:coll_problem}.
Although these two methods are mathematically equivalent, C-RONS proves to be numerically more stable.

For the Monte Carlo approximation, one draws a random sample $\{\vc x_k\}_{k = 1}^{N}$ from the spatial domain $D$ to approximate the inner products which appear in equations~\eqref{eq:metricT} and~\eqref{eq:rhs}. More specifically, taking the Hilbert space $H$ to be the space of square-integrable functions $L^2(D)$, the Monte Carlo approximations for the metric tensor $M$ and right-hand side vector $\vc f$ are given by
\begin{align}
M_{ij}(\vc q) \approx \bar M_{ij}(\vc q)  := \frac{|D|}{N} \sum_{k = 1}^N \frac{\partial \hat u}{\partial q_i}(\vc x_k, \vc q) \frac{\partial \hat u}{\partial q_j}(\vc x_k, \vc q),& \quad &i,j \in \{1,...,n\}, \nonumber\\
f_{i}(\vc q) \approx \bar f_{i}(\vc q)  :=  \frac{|D|}{N} \sum_{k= 1}^N \frac{\partial \hat u}{\partial q_i}(\vc x_k,\vc q) F(\hat u) \bigg|_{\vc x = \vc x_k},& \quad &i\in \{1,...,n\},
\label{eq:qdot_unconst_MC}
\end{align}
where $\bar M\in\mathbb R^{n\times n}$ and $\bar{\vc f}\in\mathbb R^n$ denote the Monte Carlo approximations and $|D|$ denotes the size of the spatial domain, assuming that it is bounded. 
Du and Zaki~\cite{du21} draw the samples $\{\vc x_k\}_{k = 1}^{N}$ from a uniform distribution. Bruna et al.~\cite{bruna22} showed that drawing the samples from an adaptive distribution that depends on the approximate solution $\hat u$ may lead to more accurate solutions. In either case, the unconstrained RONS~\eqref{eq:qdot_unconst} with Monte Carlo approximation can be written as 
\begin{equation}
	\bar M(\vc q) \dot{\vc q} = \bar{\vc f}(\vc q).
	\label{eq:RONS_MC}
\end{equation}

%\begin{equation}
%\bigg( \int_{\R^n} [\grad_{\vc q} \hat u] [\grad_{\vc q} \hat u]^\top  \ \id \vc x \bigg) \dot{\vc q}= \int_{\R^n} [\grad_{\vc q} \hat u] F(\hat u)\ \id \vc x.
%\label{eq:qdot_unconst_L2}
%\end{equation}

Equation~\eqref{eq:qdot_unconst_MC} reveals a close relation between C-RONS and the Monte Carlo approximation of RONS. Note that the Monte Carlo approximation of the metric tensor satisfies $\bar M = (|D|/N)\tilde{M}^\top  \tilde{M}$, where $\tilde{M}$ is the collocation matrix~\eqref{eq:M_coll}. Similarly, the Monte Carlo approximation of the right-hand side vector field $\vc f$ satisfies $\bar{\vc f} = (|D|/N)\tilde{M}^\top \tilde{\vc f}$,  where $\tilde{\vc f}$ is the C-RONS vector field~\eqref{eq:f_coll}.
Therefore, the Monte Carlo approximation of RONS~\eqref{eq:RONS_MC} can be equivalently written as
\begin{equation}
	\tilde{M}^{\top} \tilde{M} \dot{\vc q} =  \tilde{M}^{\top} \tilde{\vc f},
	\label{eq:MC_CRONS}
\end{equation}
where $\tilde M$ is the collocation matrix.

The following theorem shows that, if the C-RONS matrix $\tilde{M}$ has full column rank, then the Monte Carlo approximation of RONS and the C-RONS equation~\eqref{eq:coll_solution} are mathematically equivalent.
\begin{thm} \label{thm:coll_MC}
	If the C-RONS matrix $\tilde{M}$ has full column rank, then the Monte Carlo approximation of RONS~\eqref{eq:RONS_MC} is equivalent to C-RONS equation~\eqref{eq:coll_solution}.
\end{thm}
\begin{proof}
	First recall that the Monte Carlo approximation~\eqref{eq:RONS_MC} is equivalent to equation~\eqref{eq:MC_CRONS}. If $\tilde{M}$ is full column rank, then $\tilde{M}^\top \tilde{M}$ is invertible and therefore we have $\dot{\vc q} = ( \tilde{M}^{\top} \tilde{M} )^{-1} \tilde{M}^{\top}\tilde{\vc f}(\vc q)$.
	On the other hand, since $\tilde{M}$ has full column rank, the psuedoinverse of $\tilde{M} $ is given explicitly by $\tilde{M}^{+} = ( \tilde{M}^{\top} \tilde{M} )^{-1} \tilde{M}^{\top}$. Substituting this expression in C-RONS equation~\eqref{eq:coll_solution}, we conclude that the Monte Carlo approximation of RONS and C-RONS lead to the same equation for $\dot{\vc q}$.
\end{proof}

\begin{rem}\label{rem:condNum}
We note an important distinction between the exact result of Theorem~\ref{thm:coll_MC} and its numerical implementation.
Although Theorem~\ref{thm:coll_MC} states that C-RONS equation~\eqref{eq:rons_colloc} and the Monte Carlo approximation of RONS~\eqref{eq:MC_CRONS} are equivalent in exact arithmetic, the C-RONS equation is numerically better conditioned than the Monte Carlo approach. To see this, note that $\kappa(	\tilde{M}^{\top} \tilde{M} ) = [\kappa(\tilde{M} )]^2$, where the condition number is defined as $\kappa(\tilde{M}) := \sigma_{\text{max}} (\tilde{M}) / \sigma_{\text{min}} (\tilde{M})$ with $\sigma_{\text{max}}$ and $\sigma_{\text{min}}$ denoting the maximal and minimal nonzero singular values of the matrix, respectively.
Thus, it is numerically more stable to solve the C-RONS equation to obtain $\dot{\vc q}$. In other words, although $\tilde M^+ =  (\tilde{M}^{\top} \tilde{M} )^{-1} \tilde{M}^{\top}$ in exact arithmetic, it is well-known that this formula is numerically sensitive. Instead, we use the singular value decomposition of $\tilde M$ to compute its pseudoinverse which is numerically more stable~\cite{Ipsen2009}. 
In contrast, using the Monte Carlo approximation~\eqref{eq:MC_CRONS}, one must inevitably work with the matrix $\bar M \propto \tilde M^\top\tilde M$ which in practice tends to have a significantly larger condition number than $\tilde M$.
\end{rem}

Theorem~\ref{thm:coll_MC} sheds light on the unreasonable effectiveness of the Monte Carlo approximation applied to RONS. 
As we show in section~\ref{sec:KS} below, using Monte Carlo integration, with a relatively small samples size $N=128$, captures the behavior of the Kuramoto--Sivashinsky PDE reasonably well. This is surprising because the sample size is too small to accurately approximate the integrals involved in the metric tensor $M$ or the right-hand side vector $\vc f$ (see section~\ref{sec:KS} for a quantitative comparison). Yet, the approximate solution is reasonably close to a true solution of the Kuramoto--Sivashinsky equation. Theorem~\ref{thm:coll_MC} shows that this accuracy is not owed to the accuracy of the Monte Carlo approximation; rather it is due to the fact that this approximation, although disguised as a Monte Carlo method, is in fact a collocation method. As such, it minimizes the approximation error at $N=128$ collocation points. We discuss this point in greater detail in section~\ref{sec:KS}.

\subsection{Regularized RONS}
 \label{sec:tikh}
In sections \ref{sec:sym_comp} and \ref{sec:colloc}, we developed two efficient methods to construct the RONS equations.
The next step is to solve the resulting ODEs in order to evolve the parameters $\vc q(t)$ of the shape-morphing approximation $\hat u(\vc x,\vc q(t))$.
In our experience, when the number of parameters $n$ is large, the RONS equations may become stiff. As a result, using explicit schemes for numerical integration leads to exceedingly small time steps. 

To overcome this problem, Refs~\cite{du21,bruna22} use implicit time integration which unfortunately introduces a different set of issues. Namely, implicit methods require solving a nonlinear system at every time step which adds to the computational cost of RONS. Furthermore, the iterative methods for solving the nonlinear system are not guaranteed to converge~\cite{Keller1992,faraz_adjoint}. In order to address the possible stiffness of RONS equations, while avoiding implicit schemes, we introduce a regularized version of RONS.

\begin{comment}
\subsection{Truncated Singular Value Decomposition}

When computing the psuedoinverse of the matrix $\tilde{M}$, we first calculate the singular value decomposition of the metric tensor so that $\tilde{M} = USV^\top$.
The psuedoinverse is then calculated by $\tilde{M}^+ = VS^{+}U^\top$, where $S^{+}$ is computed by taking the reciprocal of each nonzero entry of $S$ and leaving the zeros in place.
Thus, small singular values become large entries in $S^{+}$ and cause to numerical instability.
Additionally, numerical error whenn calculating the SVD can cause singular values which should be zero to be small.

One approach to avoid this instability is to simply truncate small singular values when computing $M^+$.
This helps to increase both the speed and stability when calculating the psuedoinverse.
Truncating too many singular values will significantly alter the dynamics of the reduced-order equations, and so the tolerance for truncating small singular values must be chosen carefully.
Regardless, we obtain accurate and fast simulations with this criteria.
\end{comment}

More specifically, we add a Tikhonov penalization term~\cite{calvetti2000,golub1999} to the cost function~\eqref{eq:costfunctional}, and define the regularized cost function,
\begin{equation}
	\hat{  \mathcal J} ( \vc q, \dot{ \vc q }) = \frac{1}{2} \|\hat u_t-F(\hat u)\|_H^2 + \frac{1}{2} \| \Gamma \dot{ \vc q } \|_2^2,
	\label{eq:costfunctional_tikh}
\end{equation}
where the full-rank Tikhonov matrix $\Gamma \in \R^{P \times n}$ ($P\geq n$) is to be specified.
We then consider the constrained minimization problem,
\begin{align}\label{eq:min_cons_tikh}
	& \min_{\dot{\vc q} \in \R^n} \  \hat{  \mathcal J} ( \vc q, \dot{ \vc q }) \nonumber\\
	&\mbox{subject to}\quad I_{ k }(\vc q (t))= I_{ k }(\vc q (0)) , \quad k =1,2,...,m, \quad \forall t\geq 0.
\end{align}
As before, the constraints ensure that the resulting solution $\hat u(\vc x,\vc q(t))$ conserves the first integrals $I_k$.
The following theorem gives the explicit form of a minimizer to the regularized optimization problem~\eqref{eq:min_cons_tikh}.

\begin{thm} \label{thm:tikh_soln}
	If $\Gamma$ has full column rank and the constraint gradients  $\nabla I_1(\vc q), \nabla I_2(\vc q)$, $\cdots$, $\nabla I_m(\vc q)$ are linearly independent, then the solution to the regularized minimization problem~\eqref{eq:min_cons_tikh} satisfies
	\begin{equation}
		(M(\vc q) + \Gamma^\top \Gamma )\dot{  \vc q } = \vc f (\vc q) - \sum_{i=1}^m  \hat \lambda_k \nabla I_k(\vc q).
		\label{eq:minimizer_regularized}
	\end{equation}
	The Lagrange multipliers $\hat {\pmb \lambda} =  ( \hat\lambda_1, ... , \hat\lambda_m )^\top$ are determined through the linear system
	\begin{equation}
		\hat C(\vc q)  \hat{\pmb \lambda} = \hat{ \vc b }(\vc q) ,
		\label{eq:lambda_eq_tikh}
	\end{equation}
	where $\hat C$ is the \textit{regularized constraint matrix} with entries,
	\begin{equation}
		\hat C_{ij} = \langle \grad I_j , (M + \Gamma^\top \Gamma )^{-1} \grad I_i\rangle,\quad i,j\in\{1,2,\cdots,m\}, 
		\label{eq:constMat_tikh}
	\end{equation}
	and the vector $\hat {\vc b} = (\hat b_1, \hat b_2,\cdots, \hat b_m)^\top\in\R^m$ is given by
	\begin{equation}
		\hat b_i = \langle \grad I_i ,  (M + \Gamma^\top \Gamma )^{-1} \vc f \rangle,\quad i=1,2,\cdots, m.
		\label{eq:b}
	\end{equation}
\end{thm}
\begin{proof}
	See appendix \ref{sec:tikh_proof}.
\end{proof}

We refer to equation~\eqref{eq:minimizer_regularized} as the \emph{regularized RONS} equation.

\begin{rem}
	In addition to alleviating the stiffness of the RONS equations, the regularization also relaxes the assumptions needed on the shape-morphing solution $\hat u$. Note that the metric tensor $M$ is symmetric by definition. It is also positive semi-definite because, for all $\pmb \xi \in \R^n$, we have
	\begin{equation}
		\langle  \pmb \xi , M \pmb \xi \rangle = 	\sum_{i = 1 }^n \sum_{j= 1 }^n  \bigg \langle  \frac{\partial \hat u}{ \partial  q_i } \xi_i  , \frac{\partial \hat u}{ \partial  q_j } \xi_j \bigg \rangle_H = \bigg \| \sum_{i = 1 }^n\frac{\partial \hat u}{ \partial  q_i } \xi_i \bigg \|_H^2 \geq 0.
	\end{equation}
	In Ref.~\cite{anderson2021}, to ensure that $M$ was positive definite and therefore invertible, we required the assumption that the approximate solution was an immersion (see Lemma 1 of~\cite{anderson2021}), i.e.,
	\begin{equation}
		\dim \bigg( \text{span} \bigg \{ \frac{\partial \hat u}{ \partial  q_1 } , \frac{\partial \hat u}{ \partial  q_2 }, ... , \frac{\partial \hat u}{ \partial  q_n } \bigg \} \bigg) = n.
	\end{equation}
	In regularized RONS, we do not require the immersion assumption. Note that, for regularized RONS, we only need the invertibility of $M + \Gamma^\top \Gamma$ which is always guaranteed.
	This is because $M$ is positive semi-definite, and $\Gamma^\top\Gamma$ is symmetric positive definite. Therefore, $M + \Gamma^\top \Gamma$ is symmetric positive definite and invertible, regardless of whether the approximate solution $\hat u$ is an immersion.
\end{rem}

We can similarly apply Tikhonov regularization to C-RONS. Recall the least squares problem~\eqref{eq:coll_problem} which arises for the collocation point method, and consider its regularized counterpart,
\begin{align}\label{eq:min_cons_coll_tikh}
	& \min_{\dot{\vc q} \in \R^n} \ \|  \tilde{M} \dot{\vc q}  - \tilde{ \vc f} \|_2^2 + \| \Gamma \dot{\vc q}\|_2^2 \nonumber\\
	&\mbox{subject to}\quad I_{ k }(\vc q (t))= I_{ k }(\vc q (0)) , \quad k =1,2,...,m, \quad \forall t\geq 0.
\end{align}
The following theorem gives an explicit expression to the solution of the this optimization problem.

\begin{thm}
If $\Gamma$ has full column rank and the constraint gradients  $\nabla I_1(\vc q), \nabla I_2(\vc q)$, $\cdots$, $\nabla I_m(\vc q)$ are linearly independent, the minimizer to the constrained optimization problem~\eqref{eq:min_cons_coll_tikh} satisfies
 \begin{equation}
 	(\tilde M^\top \tilde M + \Gamma^\top \Gamma )\dot{  \vc q } = \tilde M^\top  \tilde{\vc f} - \sum_{i=1}^m  \tilde \lambda_k \nabla I_k(\vc q),
 \end{equation}
which we refer to as the regularized C-RONS equation.
The Lagrange multipliers $\tilde {\pmb \lambda} =  ( \tilde\lambda_1, ... , \tilde\lambda_m )^\top$ are determined through the linear system
 \begin{equation}
 	\tilde C(\vc q)  \tilde{\pmb \lambda} = \tilde{ \vc b }(\vc q) ,
 \end{equation}
 where the regularized constraint matrix $\tilde C$ has entries,
 \begin{equation}
 	\tilde C_{ij} = \langle \grad I_j , (\tilde M^\top \tilde M  + \Gamma^\top \Gamma )^{-1} \grad I_i\rangle,\quad i,j\in\{1,2,\cdots,m\}, 
 	\label{eq:constMat_tikh}
 \end{equation}
 and the vector $\tilde {\vc b} = (\tilde b_1, \tilde b_2,\cdots, \tilde b_m)^\top\in\R^m$ is given by
 \begin{equation}
 	\tilde b_i = \langle \grad I_i ,  (\tilde M^\top \tilde M + \Gamma^\top \Gamma )^{-1} \tilde M^\top  \tilde {\vc f} \rangle,\quad i=1,2,\cdots, m.
 \end{equation}
\end{thm}
\begin{proof}
	The proof of this theorem is very similar to that of Theorem~\ref{thm:tikh_soln} and therefore is omitted here for brevity.
\end{proof}
 
As before, invertibility of $(\tilde M^\top \tilde M + \Gamma^\top \Gamma )$ is guaranteed by the fact that $\tilde M^\top \tilde M$ is symmetric, positive semi-definite and $\Gamma^\top \Gamma $ is positive definite. In the numerical  examples presented in section~\ref{sec:num_ex}, we take the matrix $\Gamma$ to be a multiple of the identity matrix so that $\Gamma^\top \Gamma = \alpha I$, where $\alpha>0$ is a prescribed regularization parameter.

\section{Numerical results}
\label{sec:num_ex}
In this section we present two numerical examples: the Fokker--Planck equation and the Kuramoto--Sivashinsky equation. The Fokker--Planck equation demonstrate the computational advantages of using symbolic RONS as introduced in section~\ref{sec:sym_comp}. The Kuramoto--Sivashinsky equation demonstrates the benefits of using the collocation point method (section~\ref{sec:colloc}) over Monte Carlo integration. In both numerical examples, we also discuss the advantages of regularization as described in section~\ref{sec:tikh}.
We carried out our computations on a 2019 Macbook Pro with a 1.7 GHz Quad-Core Intel Core i7 processor.
Time integration for the Kuramoto--Sivashinsky equation was carried out using an explicit adaptive Runge-Kutta scheme, i.e., Matlab's \texttt{ode45}~\cite{dormand1980}. For the Fokker--Planck equation, we used an explicit adaptive multi-step solver, i.e., Matlab's \texttt{ode113}~\cite{shampine1997}. 

\subsection{Fokker--Planck equation}
In this section we consider the Fokker--Planck equation in eight dimensions.
We demonstrate that applying RONS with symbolic computation provides solutions which are several orders of magnitude more accurate and faster than the adaptive Monte Carlo sampling technique used in~\cite{bruna22}.
We also demonstrate the importance of enforcing conserved quantities to obtain accurate approximate solutions of the Fokker--Planck equation.

Following Bruna et al.~\cite{bruna22}, we consider $d$ interacting particles whose motion is governed by the the system of stochastic differential equations (SDEs),
\begin{equation}
\id X_i = g(t, X_i) \id t +  \sum_{ j = 1 }^d K(X_i, X_j) \id t +  \sqrt{2\nu} \id W_i, \quad i = 1,2,...,d.
\label{eq:SDE}
\end{equation}
%\dot{X}_i  &= a(t) - X_i + \frac{\alpha}{d}  \sum_{ j = 1 }^d (X_j - X_i) + \nu \dot{W}.
Here $X_i(t)$ denotes the position of the $i$-th particle at time $t$, $g:[0,\infty) \times  \R \rightarrow \R $ is a forcing term, $K:\R \times  \R \rightarrow \R$ describes the pairwise interactions between particles, $\nu$ is a positive diffusion constant, and $W_i$ is a standard Wiener process.
The Fokker--Planck equation is a deterministic PDE which describes the evolution of the probability density function (PDF), $p(\vc x, t)$, for the location of the particles.
The Fokker--Planck equation associated with~\eqref{eq:SDE} reads
\begin{equation}
\frac{\partial p}{ \partial t } = \sum_{ i = 1 }^{d}  - \frac{\partial  }{\partial x_i} \bigg[ \bigg( g(t,x_i)+ \sum_{ j = 1 }^{d} K(x_i,x_j) \bigg) p \bigg]  + \nu \frac{\partial^2 p }{\partial x_i^2}.
\label{eq:FP}
\end{equation}

As the spatial dimension $d$ grows, solving the Fokker--Planck equation using conventional discretization methods becomes prohibitive~\cite{Spiliopoulos2018,Weinan2018}. Alternatively, one may seek to approximate the density $p(\vc x,t)$ using Monte Carlo simulations of the original SDE~\eqref{eq:SDE}. This also becomes prohibitively expensive in higher dimensions since exceedingly large samples are required.
Here, we use RONS to directly approximates $p(\vc x, t)$, bypassing the need for any Monte Carlo simulations of the SDE or spatial discretization of the PDE.

As in \cite{bruna22}, we choose the functions $g$ and $K$ to be 
\begin{equation}
g(t, x_i) = a(t) - x_i, \quad K(x_i,x_j) = \frac{\alpha}{d} (x_j - x_i),
\label{eq:FP_gandk}
\end{equation}
which correspond to particles in a harmonic trap centered in each spatial coordinate at $a(t)$ while the particles also attract each other.  
A significant advantage of this choice is that we can obtain analytical expressions for the mean and covariance of each particle to use as a benchmark for our approximate solutions.
Taking the expected value of the SDE~\eqref{eq:SDE} with our choices of $g$ and $K$ given in~\eqref{eq:FP_gandk}, we obtain the following expressions for the mean $\bar{X}_i = \mathbb{E}[X_i]$ of each particle
\begin{equation}
\dot{\bar{X}}_i = a(t) - \bar{X}_i + \frac{\alpha}{d}  \sum_{ j = 1 }^d (\bar{X}_j - \bar{X}_i),\quad i = 1,2,...,d.
\label{eq:trap_mean}
\end{equation}
Similarly, we have the following expressions for entries of the matrix $\Sigma_{ij} = \mathbb{E}[X_i X_j]$,
\begin{align}
\dot{\Sigma}_{ij} &= a(t) ( \bar{X}_j + \bar{X}_i ) - 2(1+\alpha)\Sigma_{ij} + \frac{\alpha}{d} \sum_{ l = 1 }^d  ( \Sigma_{lj} + \Sigma_{li} ) + 2  \nu \delta_{ij},& \quad &i,j \in \{1,2,...,d\} ,
\label{eq:trap_sigma}
\end{align}
where  $\delta_{ij}$ denotes the Kronecker delta. The covariance matrix $\Sigma_{ij}-\bar X_i\bar X_j$ can then be computed using the solutions to equations~\eqref{eq:trap_mean} and~\eqref{eq:trap_sigma}.

We consider the same initial condition and parameter values as in Ref.~\cite{bruna22}. More specifically, we take the Gaussian initial condition,
\begin{equation}
p(\vc x, 0) = (2 \pi )^{-d/2}\det(\Sigma)^{-1/2} \exp \bigg[ -\frac{1}{2} (\vc x - \pmb \mu )^\top \Sigma^{-1} (\vc x - \pmb \mu  ) \bigg],
\end{equation}
where the initial mean $\pmb \mu \in \R^d $ is given by $\mu_i =  0.9 + 2.1(i-1)/(d-1)$ for $i = 1,...,d,$ and the initial covariance $\Sigma \in \R^{d\times d}$ is the diagonal matrix $\Sigma = \text{diag}(0.1,0.1,...,0.1)$. 
The remaining parameters are given by $a(t) = 1.25(\sin(\pi t) + 1.5) $, $\alpha = 0.25$, $d = 8$, and $\nu = 0.01$. 

We approximate the solution of the Fokker--Planck equation~\eqref{eq:FP} using symbolic RONS as described in section \ref{sec:sym_comp}.
For the shape-morphing approximate solution~\eqref{eq:ansatz}, we choose
\begin{equation}
\hat{p} (\vc x, \vc q(t)) = \sum_{ i = 1 }^r A_i^2(t) \exp \left[ - w_i^2(t) |\vc x - \vc c_i(t)|^2 \right],
\label{eq:gauss_ansatz_FP}
\end{equation}
where the mode function $\phi$ is a Gaussian, the shape parameters are $\pmb\beta_i=(w_i,\vc c_i)\in\mathbb R^9$ and the amplitudes are $\alpha_i = A_i^2$.
We square the amplitudes to ensure that the approximate PDF $\hat p$ is non-negative.
The wights $w_i(t)$ control the standard deviation of each Gaussian since $|w_i|^{-1}$ is proportional to the standard deviation of the $i$-th Gaussian. Finally, the vector $\vc c_i(t)\in \R^d$ determines the $i$-th Gaussian's mode. Therefore, the parameters of the approximate solution are $\vc q = \{A_i,w_i,\vc c_i\}_{i=1}^r$, resulting in a total of $n=r(d+2)$ parameters.

Note that since the solution $p$ is a PDF, its integral over the entire domain $\mathbb R^d$ must be equal to one for all times. This constitutes a conserved quantity for the Fokker--Planck equation which can be easily enforced in RONS.
We ensure that the total probability of the approximate solution $\hat p$ is unity by enforcing the conserved quantity,
\begin{equation}
I_1(\vc q(t)) := \int_{\R^d} \hat p(\vc x, \vc q(t)) \, \id \vc x =  \pi^4\sum_{ i = 1 }^r \frac{A_i^2(t)}{w_i^8(t)} = 1,
\label{eq:FP_const}
\end{equation}
for all $t\geq 0$.

Assuming that the Hilbert space $H$ is the space of square integrable functions over $\R^d$, we use the symbolic RONS with the Gaussian approximate solution~\eqref{eq:gauss_ansatz_FP} to form the RONS equation.
Since our initial condition $p(\vc x,0)$ is a Gaussian and the approximate solution is a sum of Gaussians, there are infinitely many choices of parameter $\vc q(0)$ with which the approximate solution can exactly represent the initial condition.
We choose to represent the initial condition by giving all Gaussians in the approximate solution $\hat p$ the same mean and covariance as the initial condition, and then equally distributing the amplitude of the initial condition between each of the $r$ Gaussians in the approximate solution.
More explicitly, we choose initial parameter values $A_i^2(0) = (2 \pi \times 0.1 )^{-4}r^{-1}$, $ w_i^2(0) =  (2\times 0.1 )^{-1}$, and $\vc c_i(0) = \pmb \mu$.
With this choice of initial parameters, the metric tensor $M$ is not invertible at the initial time and so we use the Moore-Penrose pseudoinverse $M^+$ when solving the RONS equation~\eqref{eq:qdot_const}.

\subsubsection{Symbolic RONS without regularization}
\begin{figure}
	\centering
	\includegraphics[width=0.49\textwidth]{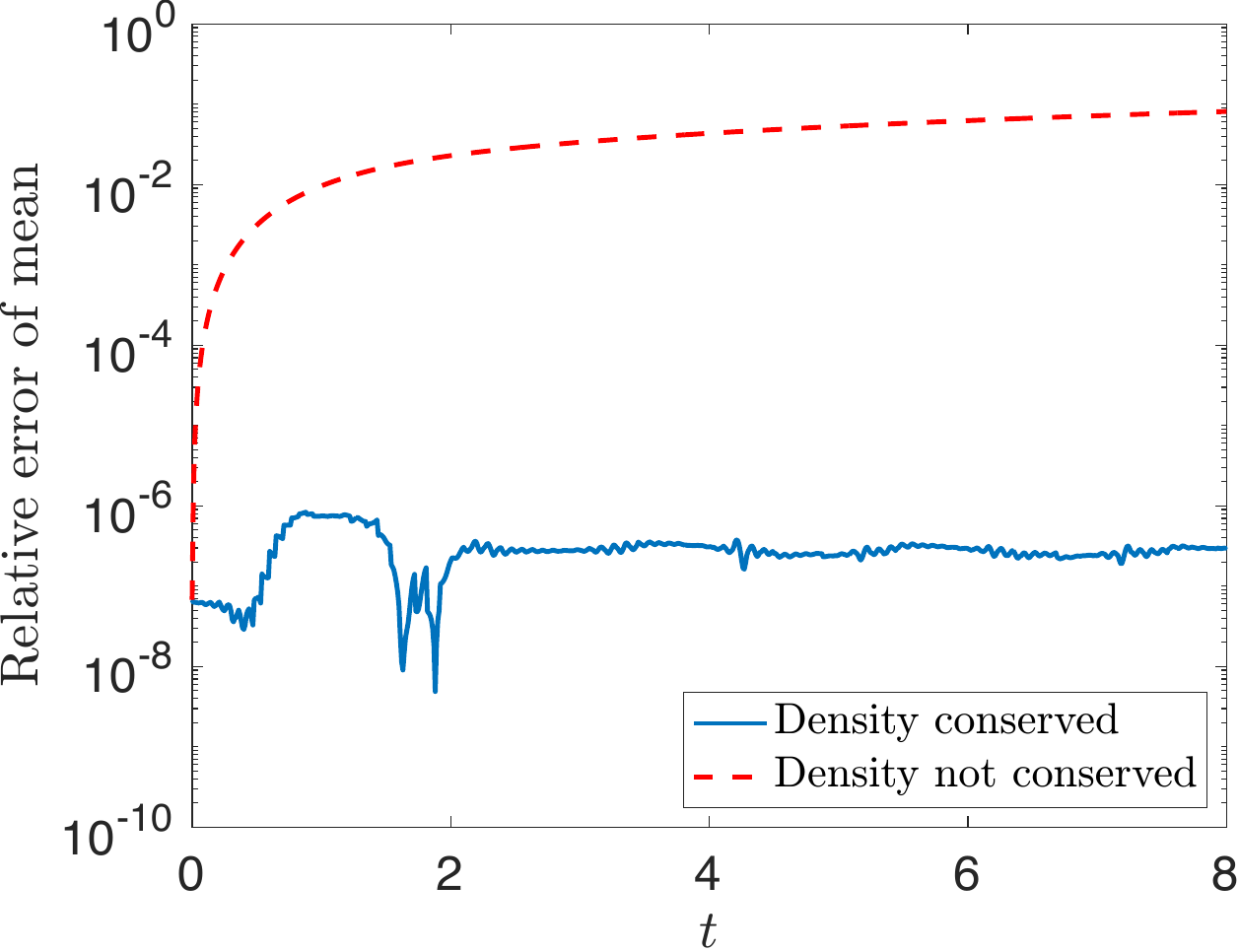}
	\includegraphics[width=0.49\textwidth]{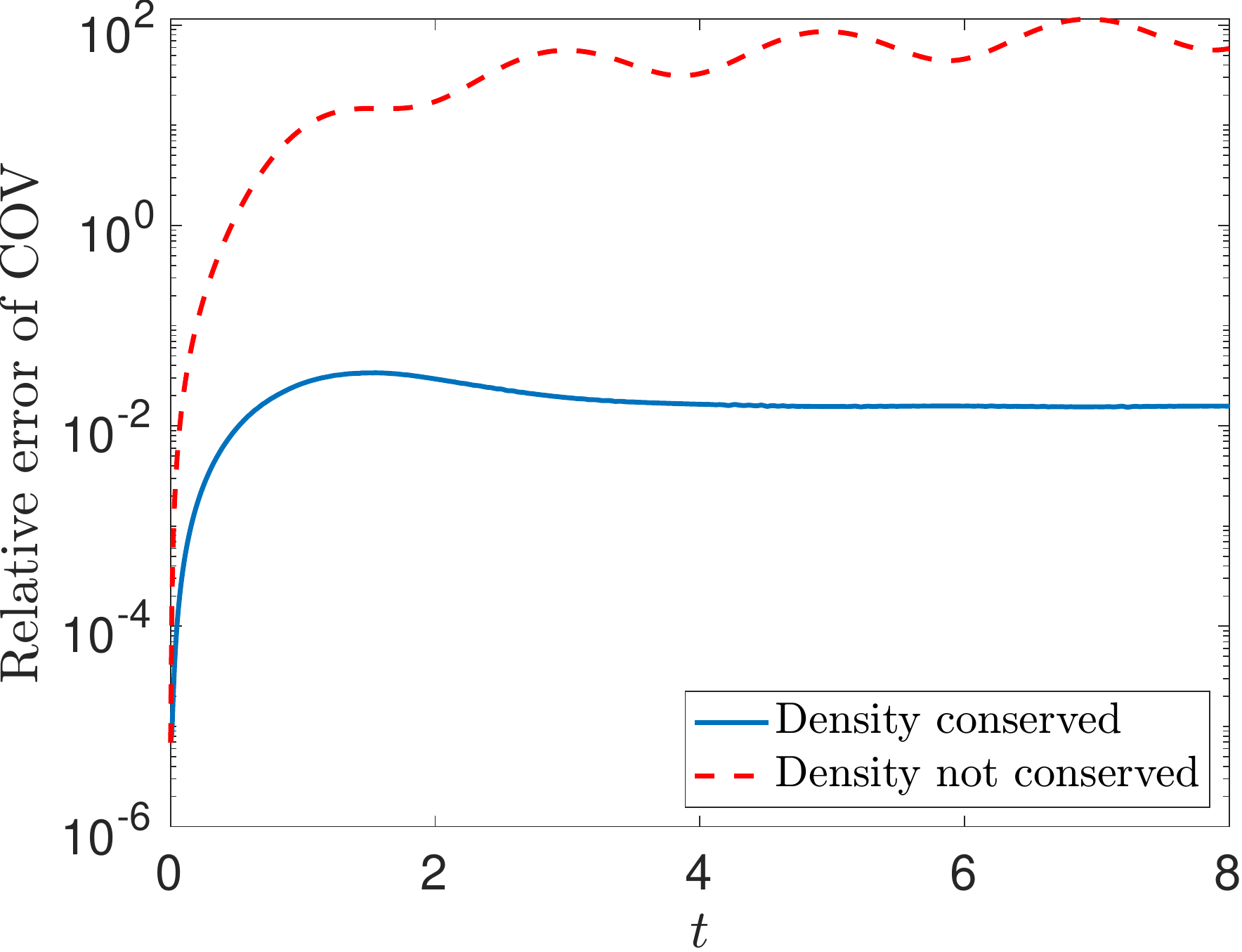}
	\caption{Relative error for mean and covariance of the S-RONS solution applied to the Fokker--Planck equation for the harmonic trap using two Gaussians in the approximate solution $(r=2)$.}
	\label{fig:Harmonic_2Gauss}
\end{figure}

First, we study solutions to RONS with only two modes $(r = 2)$ in the Gaussian approximate solution~\eqref{eq:gauss_ansatz_FP}.
In this case, the resulting ODEs are not stiff and therefore regularization is not necessary. We first demonstrate the importance of enforcing conserved quantities in the reduced-order model. Figure~\ref{fig:Harmonic_2Gauss} shows the results both with and without enforcing that the total probability of the approximate solution remains constant; see equation~\eqref{eq:FP_const}.
When enforcing constant total probability, RONS captures the true mean with a relative error on the order of $10^{-6}$ and the covariance is captured with relative error of approximately $10^{-2}$  (solid blue curves in figure \ref{fig:Harmonic_2Gauss}).
In contrast, if the conservation of probability is not enforced, the relative error of the mean increases to about $10^{-1}$ and the  covariance error reaches $10^2$ (dashed red curves in figure \ref{fig:Harmonic_2Gauss}).
Therefore, enforcing the conserved quantity~\eqref{eq:FP_const} results in approximate solutions which are 4 to 5 orders of magnitude more accurate.

In addition to providing accurate solutions, the time integration for the RONS simulation takes only 0.31 seconds when using 2 Gaussians in the approximate solution. If we were to instead simulate many realizations of the SDE to approximate the PDF, the total computational time would be significantly higher.
 
 As mentioned earlier, the Fokker--Planck equation~\eqref{eq:FP} was also solved in~\cite{bruna22}, where they used an adaptive Monte Carlo sampling to estimate the inner products in RONS. In Table~\ref{table:harmonic_2gauss}, we compare the computational time and accuracy between the adaptive sampling approach of~\cite{bruna22} and our symbolic RONS, where both methods use two Gaussians in the approximate solution~\eqref{eq:gauss_ansatz_FP}. For the adaptive sampling method, we use the code that was made publicly available by Bruna et al.~\cite{Bruna2023_github}. This code uses a backwards Euler scheme for time integration together with stochastic gradient descent to solve the nonlinear system at every timestep, whereas our RONS simulations are integrated in time using an explicit scheme. As shown in Table~\ref{table:harmonic_2gauss}, we see that RONS returns significantly faster and more accurate solutions than the adaptive sampling approach. 
 
In particular, time integration using adaptive sampling takes approximately 189 minutes (more than 3 hours). The main computational cost comes from sampling and subsequent evaluation of the inner products, which has to be repeated at each time step. In contrast, time integration using symbolic RONS only takes 0.31 seconds since no sampling is required and the symbolic computations do not need to be repeated at every time step. There is the one-time cost of computing the integrals symbolically for RONS which takes approximately 13.7 minutes. However, the symbolic expressions for the RONS equation are obtained, we can integrate the equations from any initial condition without having to recompute the inner products. In other words, the solution from a different initial condition $p(\vc x,0)$ can be obtained in approximately $0.31$ seconds.
In contrast, adaptive Monte Carlo simulations need to be repeated for every initial condition and therefore the numerical integration would again take hours if we were to change the initial condition.

In addition to being faster, symbolic RONS is also more accurate. As shown in Table~\ref{table:harmonic_2gauss}, the mean of the solution is computed four orders of magnitude more accurately when using symbolic RONS compared to adaptive sampling. Moreover, the estimated covariance is two orders of magnitude more accurate when using symbolic RONS. The higher accuracy of symbolic RONS is not surprising since the inner products are computed exactly, whereas relatively large errors are accrued when adaptive Monte Carlo sampling is used.

\begin{table}[]
	\centering
		\caption{Comparison of computational time and accuracy for harmonic trap example using the adaptive sampling approach \cite{bruna22} and S-RONS. Two Gaussians are used in approximate solution for each simulation. Note that the symbolic computation only needs to be performed once; changing the initial condition or Fokker--Planck parameters does not require additional symbolic computation.}
	\begin{tabular}{|l||cccc|}
		\hline
		\textbf{Harmonic Trap $(r = 2)$} &
		\begin{tabular}[c]{@{}c@{}}Symbolic \\ computation\end{tabular} &
		\begin{tabular}[c]{@{}c@{}}Time \\ integration\end{tabular} &
		\begin{tabular}[c]{@{}c@{}}Relative error\\ of mean\end{tabular} &
		\begin{tabular}[c]{@{}c@{}}Relative error \\ of covariance\end{tabular} \\ \hline
		\multicolumn{1}{|l||}{Adaptive Sampling} &
		none &
		$ 189.1$ minutes &
		$\approx 4\times 10^{-3}$ &
		$ \approx 2$ \\ \hline
		\multicolumn{1}{|l||}{Symbolic RONS} &
		13.7 minutes &
		0.31 seconds &
		$\approx 3\times 10^{-7}$ &
		$\approx 10^{-2}$\\
		\hline
	\end{tabular}
	\label{table:harmonic_2gauss}
\end{table}

\subsubsection{Regularized symbolic RONS}
Next we consider the effect of increasing the number of modes $r$. In particular, we consider the approximate solution~\eqref{eq:gauss_ansatz_FP} with $r=30$ modes. As the number of parameters increase, the RONS equation becomes stiff. To address this issue, Bruna et al.~\cite{bruna22} use an implicit time integration scheme.
As mentioned in section~\ref{sec:tikh}, an alternative approach is to use a regularization. Here, we use regularized symbolic RONS 
with the regularization parameter $\alpha = 10^{-3}$ and compare our results to the adaptive sampling method of~\cite{bruna22} with their implicit time integrator~\cite{Bruna2023_github}.

\begin{figure}
	\centering
	\includegraphics[width=0.49\textwidth]{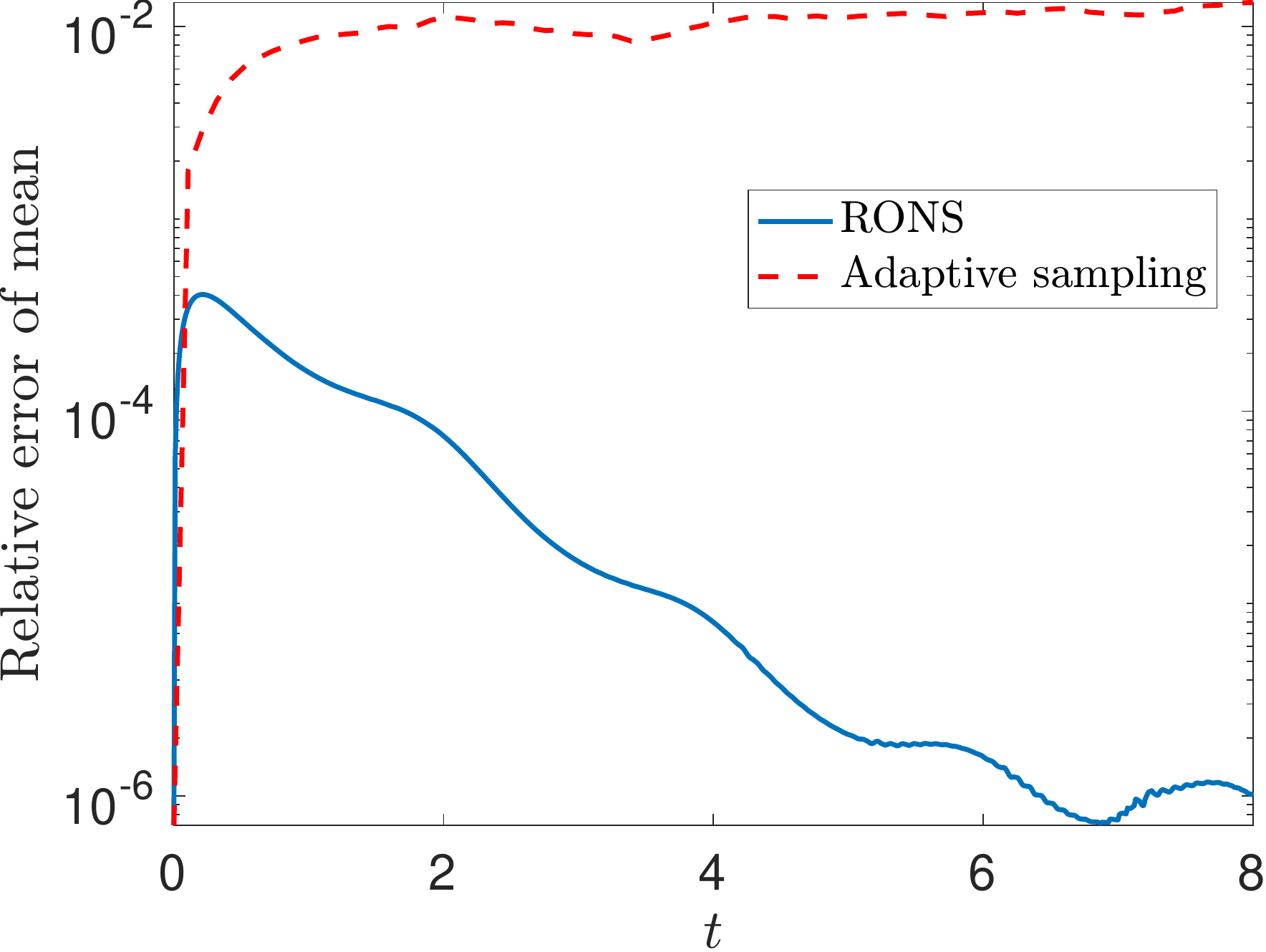}
	\includegraphics[width=0.49\textwidth]{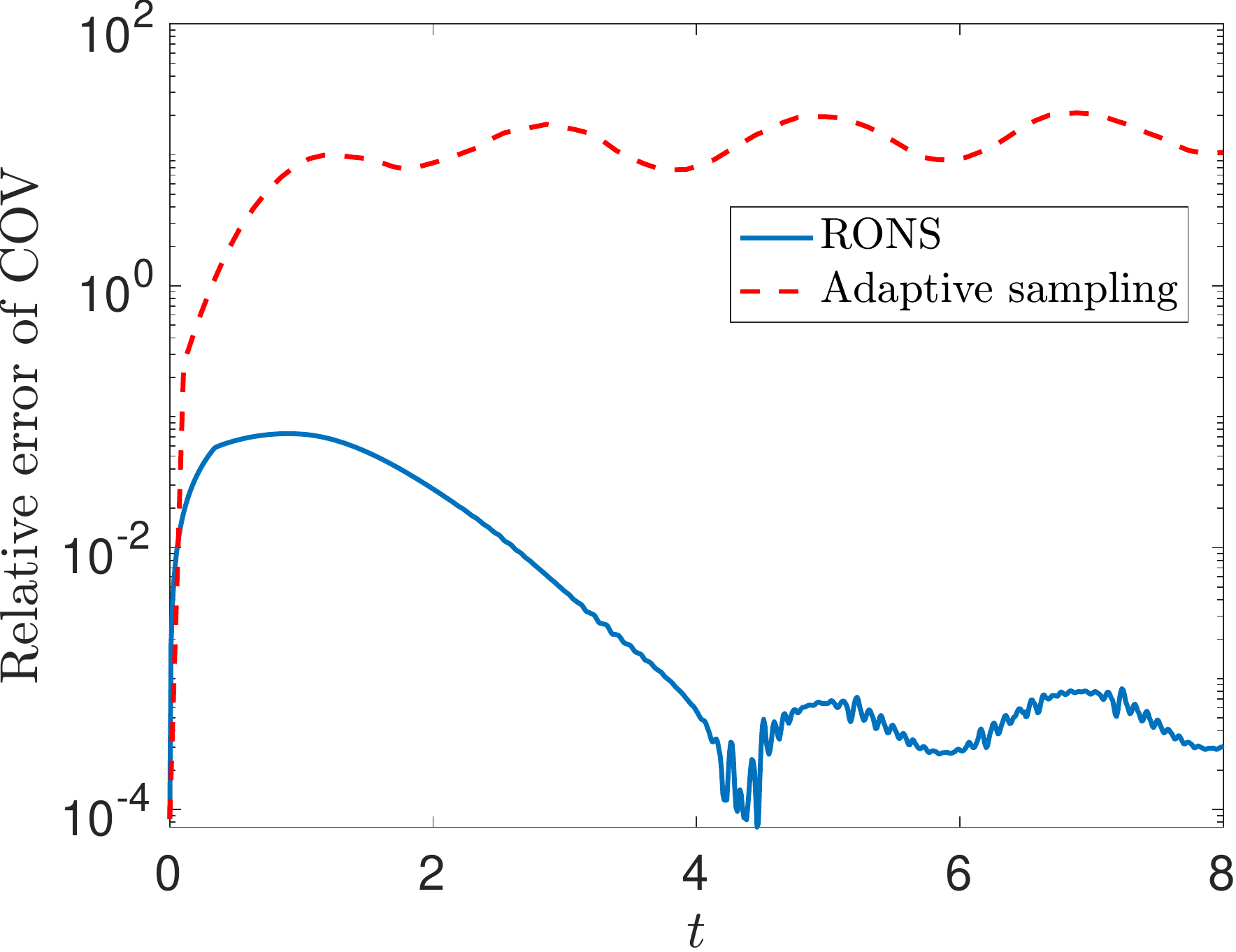}
	\caption{Relative error for mean and covariance of the S-RONS solution applied to the Fokker--Planck equation corresponding to the harmonic trap. For comparison, the same relative errors are shown for the adaptive sampling method~\cite{bruna22}.
		Thirty Gaussians are used in the approximate solution ($r=30$).
	}
	\label{fig:Harmonic_40Gauss}
\end{figure}

Figure \ref{fig:Harmonic_40Gauss} shows the relative error of the mean and covariance when applying RONS to the Fokker--Planck PDE with $r=30$ Gaussians in the approximate solution. The regularized symbolic RONS approximation matches the analytical solution well.
As the solution evolves, the relative error of the mean settles around $10^{-6}$.
We see a similar behavior in the approximate solution's covariance, where the relative error settles around $3\times 10^{-4}$ as the solution evolves.

There is a short time period around $t=0$ where the relative errors increase.
This transient increase coincides with the time needed for the particles to settle in the harmonic trap.
Initially the particles travel from their initial condition, but after a short time they settle in the harmonic trap and oscillate there. After this trapping, it becomes easier for the approximate solution to capture the true solution and therefore the relative error decreases.
Note that this initial growth was absent when using only $r=2$ modes, where no regularization was required. This demonstrates the fact that, although regularized RONS speeds us the time integration, it can lead to a deterioration of the accuracy. Nonetheless, the error is relatively small. In fact, as shown in figure~\ref{fig:Harmonic_40Gauss}, regularized symbolic RONS is about three orders of magnitude more accurate than adaptive sampling with implicit time integration.

In Table \ref{table:harmonic}, we compare the computational time and accuracy of the adaptive sampling approach of~\cite{bruna22} and regularized symbolic RONS. As in the case of 2 Gaussians, the 30-mode approximation using symbolic RONS significantly outperforms the results from adaptive sampling in both computational speed and accuracy. In particular, time integration using symbolic RONS takes slightly over one hour, whereas adaptive sampling takes over 24 hours and yet yields lower accuracy. The high computational cost of adaptive sampling is attributed to the fact that RONS inner products must be reevaluated at every time step. Furthermore, since an implicit time integration scheme is used, a nonlinear equation needs to be solved at every time step which adds to the computational cost. In contrast, regularized RONS uses an explicit scheme which does not require solving a nonlinear equation.

We emphasize that the computational cost of symbolic integration is independent of the number of modes $r$ when using symbolic RONS as described in section~\ref{sec:sym_comp}. More specifically, the symbolic expressions from $r=2$ modes can be used to evaluate all RONS terms when $r=30$, without requiring additional symbolic computation, thus the zero symbolic computational time reported in Table~\ref{table:harmonic}.

We conclude this section by remarking that S-RONS is computationally feasible for this problem because of the method developed in section~\ref{sec:sym_comp}. The Gaussian approximate solution~\eqref{eq:gauss_ansatz_FP} has $r=30$ modes with $K=10$ parameters in each mode, so that a brute force approach to RONS would have required symbolic computation of $300(300+3)/2 = 45,450$ integrals. In contrast, S-RONS requires symbolic computation of only $10(10+3)/2 = 65$ integrals (see Theorem~\ref{thm:comptrick}).
\begin{table}[]
	\centering
		\caption{Comparison of computational speed and accuracy for the harmonic trap with $r=30$ modes. We compare the adaptive sampling approach~\cite{bruna22} to our proposed method of regularized symbolic RONS.}
	\begin{tabular}{|l||cccc|}
		\hline
		\textbf{Harmonic Trap $(r = 30)$} &
		\begin{tabular}[c]{@{}c@{}}Symbolic \\ computation\end{tabular} &
		\begin{tabular}[c]{@{}c@{}}Time \\ integration\end{tabular} &
		\begin{tabular}[c]{@{}c@{}}Relative error\\ of mean\end{tabular} &
		\begin{tabular}[c]{@{}c@{}}Relative error \\ of covariance\end{tabular} \\ \hline
		\multicolumn{1}{|l||}{Adaptive Sampling} &
		none &
		$ 24.4$ hours &
		$\approx 10^{-2}$ &
		$\approx 10 $ \\
		\hline
		\multicolumn{1}{|l||}{Symbolic RONS} &
		0 minutes &
		64.2 minutes &
		$\approx 10^{-6}$ &
		$\approx 3\times 10^{-4}$\\
		\hline
	\end{tabular}
	\label{table:harmonic}
\end{table}

\subsection{Kuramoto--Sivashinsky equation} \label{sec:KS}
In this section, we consider the Kuramoto--Sivashinsky (KS) equation and approximate its solution with a shallow neural network with hyperbolic tangent activation functions. In this case, obtaining symbolic expressions for the RONS equation is not possible. Therefore, we use collocation RONS as described in section~\ref{sec:colloc}. The KS equation was also solved in~\cite{du21} using a Monte Carlo method to approximate the RONS equations. We compare our results with this Monte Carlo approach.

The Kuramoto--Sivashinsky equation is given by
\begin{equation}
\frac{\partial u }{ \partial t } = -u \frac{\partial u }{ \partial x } - \frac{\partial^2 u }{ \partial x^2 } - \frac{\partial^4 u }{ \partial x^4 }, \quad u(x,0)=u_0(x),
\label{eq:KS}
\end{equation}
where the solution $u(x,t)$ is assumed to have periodic boundary conditions over the domain $x\in[-\ell,\ell]$.
We consider the same set-up used in~\cite{du21}. In particular, we set $\ell = 10$ and consider the initial condition
\begin{equation}
	u_0(x)= -\sin \bigg(  \frac{\pi x }{\ell} \bigg), \quad x \in [-\ell,\ell].
\end{equation}
The corresponding solutions of the KS equation are known to exhibit chaotic behavior~\cite{hyman1986,Kevrekidis1990,cvitanovic2010,pathak2018}.
As the ground truth,  we use direct numerical simulations (DNS) using a Fourier pseudo-spectral method with  $2^7$ modes.

%We examine the ability of our collocation point method to approximate solutions of the Kuramoto--Sivashinsky equation.
For the KS equation, the choice of appropriate approximate solution is not as clear as in the Fokker--Planck example.
Motivated by architectures typically used in neural networks, we choose an approximate solution which is a shallow neural network with hyperbolic tangent activation function,
\begin{equation}
\hat u (x, \vc q) = \sum_{ i = 1 }^{r} A_i(t) \tanh \bigg(  w_i(t)s_i(x) + d_i(t) \bigg),
\label{eq:tanh_ansatz}
\end{equation} 
where 
\begin{equation}
s_i(x) = \sin\bigg( \frac{ \pi x }{\ell}  + c_i(t) \bigg),
\end{equation}
is a nonlinear coordinate transformation to ensure that the approximate solution $\hat u$ satisfies the periodic boundary conditions over the domain $[-\ell,\ell]$. This is a common technique discussed further in~\cite{dong21,du21,yazdani2020}. The approximate solution~\eqref{eq:tanh_ansatz} can be thought of as a neural network with a single hidden layer and $r$ nodes. Each node contains the amplitude $A_i(t)$, the weight $w_i(t)$, and the biases $c_i(t)$ and $d_i(t)$. These form the parameters of the approximate solution, $\vc q = \{A_i,w_i,c_i,d_i\}_{i=1}^r$, where the shape parameters comprise $\pmb\beta_i = (w_i,c_i,d_i)$. Typically these parameters would be obtained by a training process. However, as Du and Zaki~\cite{du21} observe, no training is required; the parameters can be evolved using the RONS equation~\eqref{eq:qdot_unconst}.

\begin{figure}
	\centering
	\includegraphics[width=\textwidth]{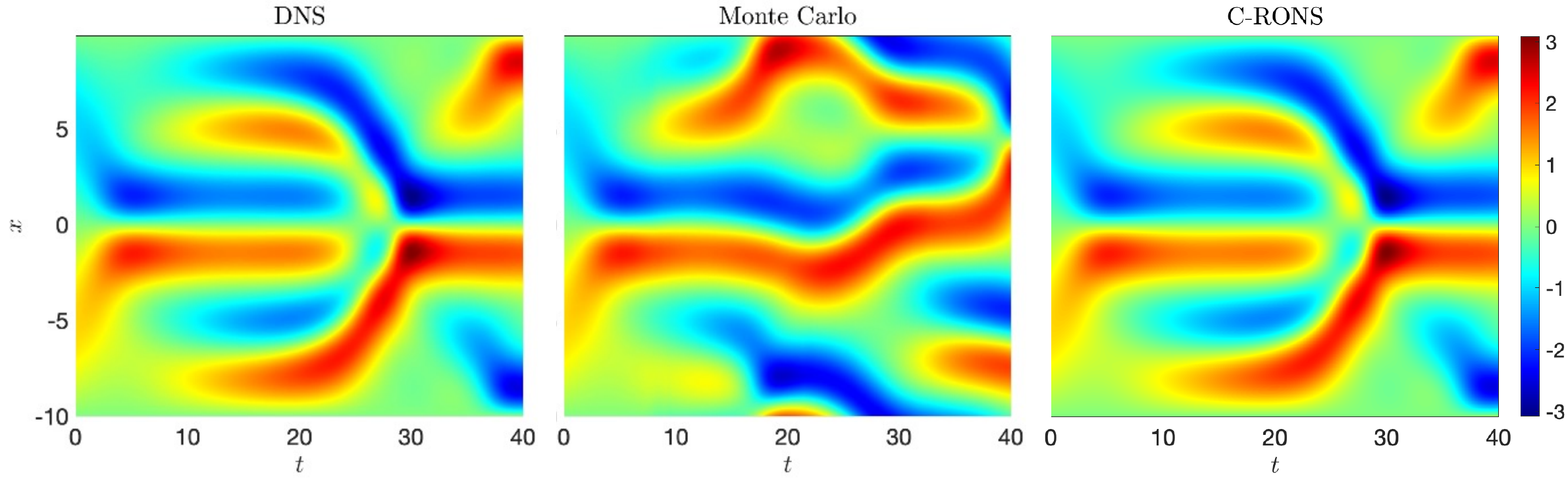}
	\caption{Simulation of Kuramoto--Sivashinsky equation using direct numerical simulation (left panel), Monte Carlo integration to approximate the RONS equation (middle panel), and regularized collocation RONS (right panel).}
	\label{fig:KS_RONS}
\end{figure}

\begin{table}[]
	\centering
		\caption{Computational time and accuracy for the Kuramoto--Sivashinsky equation using Monte Carlo integration and C-RONS.}
	\begin{tabular}{|l||lll|}
		\hline
		\textbf{Kuramoto--Sivashinsky} &
		\multicolumn{1}{c}{Monte Carlo} &
		\multicolumn{1}{c}{\begin{tabular}[c]{@{}c@{}}Regularized\\ Monte Carlo\end{tabular}} &
		\multicolumn{1}{c|}{\begin{tabular}[c]{@{}c@{}}Regularized\\ C-RONS\end{tabular}} \\ \hline
		\multicolumn{1}{|l||}{Computational time} &
		222.3 minutes &
		118.3 seconds &
		38.6 seconds\\ \hline
		\multicolumn{1}{|l||}{Relative error} &
		$1.34$ &
		$9.3\times 10^{-2}$ &
		$3.3\times 10^{-2}$\\
		\hline
	\end{tabular}
	\label{table:KS}
\end{table}

There does not exist a choice of parameters such that the sine wave initial condition $u_0$ can be exactly represented with our choice of approximate solution~\eqref{eq:tanh_ansatz}. 
To determine the initial parameter values $\vc q(0)$, we perform a least-squares fitting of the approximate solution to the initial condition, i.e., we set
\begin{equation}
\vc q(0) = \argmin_{\vc q \in \R^n} \| \hat{u}(\cdot, \vc q)  - u_0 \|_{L^2}^2.
\end{equation}
We solve this optimization problem once at the initial time to obtain the parameters $\vc q(0)$. The corresponding approximation error is less than $2 \times 10^{-7}$.

We then evolve the approximate solution~\eqref{eq:tanh_ansatz} with 10 modes $(r = 10)$ using collocation RONS method of section~\ref{sec:colloc}.
We compare our results with the Monte Carlo approximation of the inner products as proposed in~\cite{du21}. For collocation RONS we use $2^7$ equidistant collocation points. For the Monte Carlo approach, we use $2^7$ samples uniformly distributed throughout the domain $[-\ell,\ell]$. 
When applying collocation RONS, we use Tikhonov regularization as described in section~\ref{sec:tikh}, with a regularization parameter value of $\alpha = 10^{-5}$. Both methods, collocation RONS and the Monte Carlo approach, use Matlab's \texttt{ode45} for time integration.
\begin{figure}
	\centering
	\includegraphics[width=\textwidth]{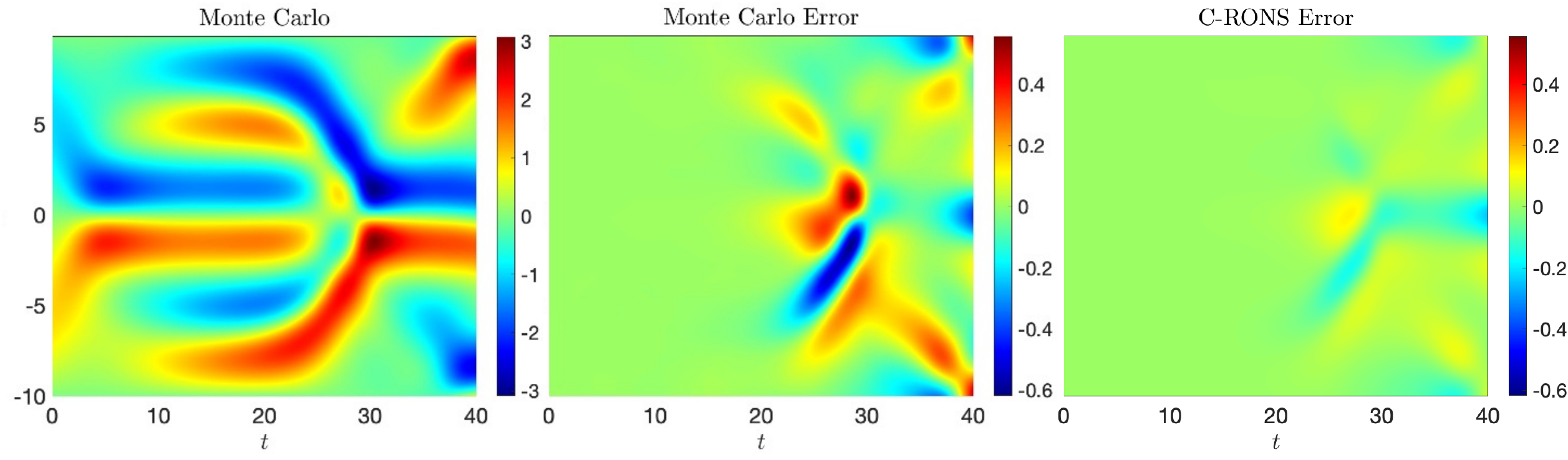}
	\caption{Simulation of Kuramoto--Sivashinsky equation using regularized Monte Carlo method (left panel). Error of the regularized Monte Carlo method (middle panel) as compared to the DNS solution. Error of regularized C-RONS is shown for comparison (right panel).}
	\label{fig:KS_MC_higherdroptol}
\end{figure}

In Figure \ref{fig:KS_RONS}, we compare the approximate solutions produced by collocation RONS and the Monte Carlo approach. 
Collocation RONS is in excellent agreement with the DNS solution, whereas the Monte Carlo approach quickly diverges from the DNS solution after approximately 10 time units. Table~\ref{table:KS} compares the computational time of these two methods. The Monte Carlo approach takes 222.3 minutes to run while the regularized collocation RONS takes only 38.6 seconds. 

The Monte Carlo method is significantly slower mainly because of the stiffness of the resulting ODEs due to the high condition number of the metric tensor $\bar M$. To demonstrate this, we also regularize the Monte Carlo approach by applying Tikhonov regularization to~\eqref{eq:RONS_MC}. As shown in Table \ref{table:KS}, regularization significantly reduces the computational time of the Monte Carlo method from 222.3 minutes to 118.3 seconds. Although this is still 3 times slower than regularized C-RONS, Tikhonov regularization greatly reduces the computational cost of the Monte Carlo method.

Interestingly, regularization also increases the accuracy of the Monte Carlo method. In Figure~\ref{fig:KS_MC_higherdroptol}, we show the error between regularized Monte Carlo method and the DNS solution. For comparison, we also show the error for regularized C-RONS.
The Monte Carlo approximation uses the regularization parameter $\alpha = 10^{-8}$ since
larger values of $\alpha$ led to numerical results which deviated significantly from the DNS.
Comparing figures~\ref{fig:KS_RONS} and~\ref{fig:KS_MC_higherdroptol}, we first note that regularization greatly improves the accuracy of the Monte Carlo approach. However, the solution obtained by regularized C-RONS still yields lower errors and only takes a third of the computational time of regularized Monte Carlo (see Table~\ref{table:KS}).
This is largely due to the poor conditioning of the matrix $\bar M$ which appears in the Monte Carlo approximation (see Remark~\ref{rem:condNum}).
\begin{figure}
	\centering
	\includegraphics[width=0.49\linewidth]{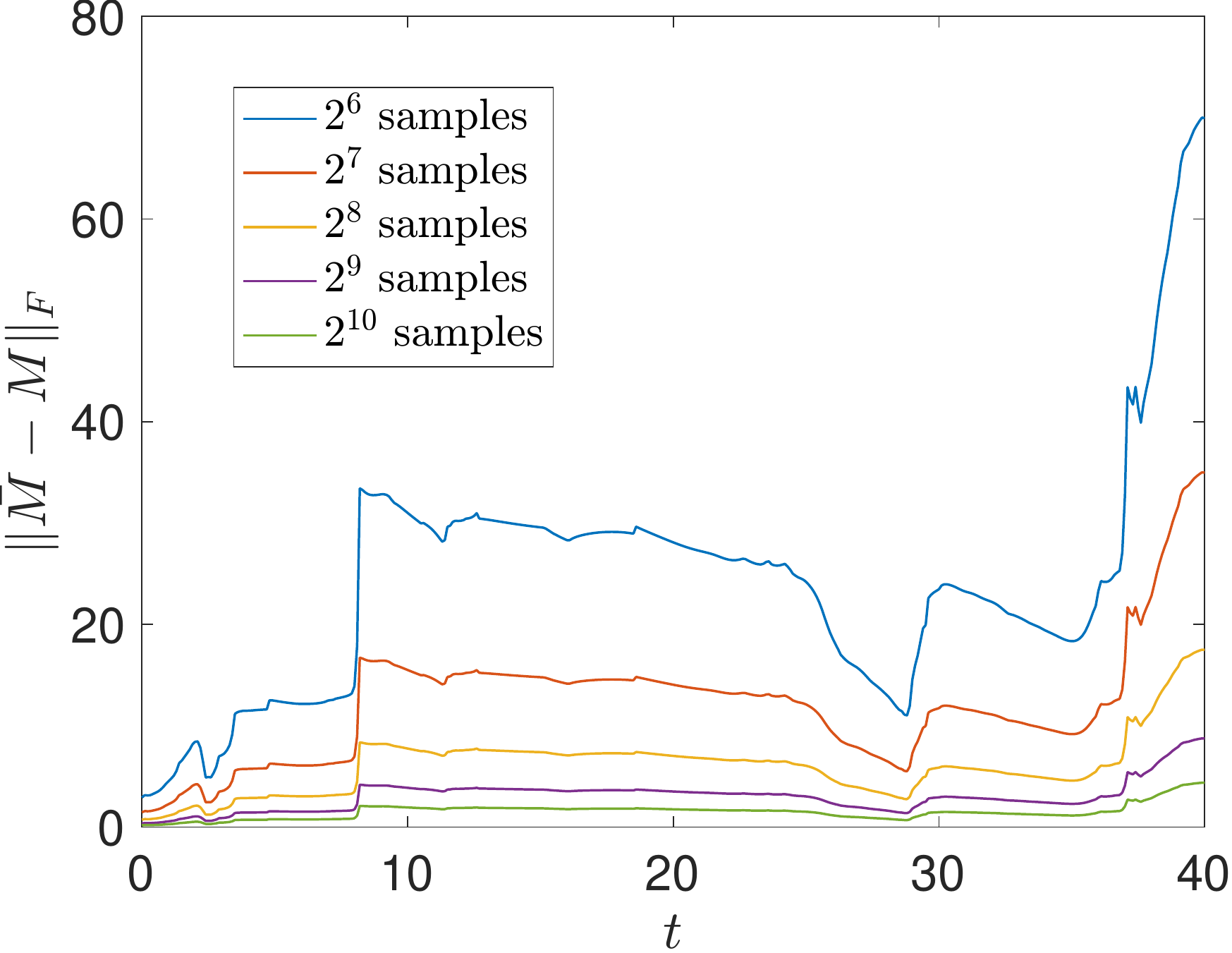}
	\caption{ Frobenius error of the Monte Carlo approximation $\bar M$ to the metric tensor for the Kuramoto--Sivashinsky equation when using increasing number of samples. Samples are all uniformly distributed on the domain $x\in [-10,10]$ and the true value of $M$ is approximated by using trapezoidal rule with $10^4$ samples.}
	\label{fig:KS_condm}
\end{figure}

The approximation errors for both methods grow over time, which is expected for a chaotic system as any error in approximating the initial condition will grow as the solution evolves. Even for short timescales, it is surprising that the Monte Carlo method is able to approximate the solution given that only $128$ samples are used. In fact, the Monte Carlo integration~\eqref{eq:qdot_unconst_MC} with 128 samples is quite inaccurate. For instance, figure~\ref{fig:KS_condm} shows the Frobenius error between the Monte Carlo approximation $\bar M$ and the true metric tensor $M$. For 128 samples, the error start around 2 and grows to approximately 18 in less than 10 time units. As a result, there is considerable error in the Monte Carlo approximation of the RONS equations when only 128 samples are used. Yet, the regularized Monte Carlo method returns a rather accurate solution as shown in figure~\ref{fig:KS_MC_higherdroptol}. Theorem~\ref{thm:coll_MC} is the key to resolving this seeming paradox. The RONS equation approximated by Monte Carlo integration is equivalent to collocation RONS.
Therefore, the Monte Carlo approximation is in fact minimizing the residual function~\eqref{eq:residual} at the sampled points.
This allows the Monte Carlo approach to obtain an accurate solution with a small sample size, despite being inaccurate as an integration method.

\section{Conclusions}\label{sec:conc}
Despite being in its infancy, RONS has already emerged as an effective method both for reduced-order modeling~\cite{anderson2021,anderson22} and for solving PDEs with neural networks without requiring any training~\cite{du21,bruna22}.
However, brute force construction of the RONS equations requires the evaluation of $\mathcal O(n^2)$ integrals, where $n$ denotes the number of time-dependent parameters in the approximate solution. Therefore, this approach becomes computationally prohibitive when a large number of parameters are required to accurately approximate the solution of the PDE. Making matters worse, the resulting ODEs tend to become stiff as the number of parameters grows. Here, we developed three methods to address these computational bottlenecks: symbolic RONS, collocation RONS, and regularized RONS.

Using symbolic computing and exploiting the structure of the RONS equations, symbolic RONS (or S-RONS, for short) drastically reduces the computational cost from $\mathcal O(n^2)$ to $\mathcal O(K^2)$ where $K\ll n$ is independent of the number of parameters $n$. Furthermore, since the equations are constructed symbolically, the S-RONS integrals do not need to be recomputed during time stepping; rather they can be evaluated by direct substitution of the updated parameters into the symbolic expressions. Applying S-RONS to the Fokker--Planck equation, we obtained 14-23 times speedup in comparison to the adaptive sampling method of~\cite{bruna22}. In addition, the accuracy of the solutions increased by several orders of magnitude.

We also developed collocation RONS (or C-RONS, for short) in case symbolic computation of the integrals are not feasible.
Rather than minimizing the error between evolution of the approximate solution and dynamics of the governing PDE over the entire spatial domain, C-RONS minimizes this error over a set of prescribed collocation points. 
Since this method does not require any symbolic computation, it is applicable to any choice of the approximate solution and any form of the PDE.
We also proved that, in exact arithmetic, C-RONS is equivalent to the Monte Carlo method proposed in~\cite{du21}. However, from a numerical standpoint, C-RONS is significantly better conditioned than the Monte Carlo approximation and thus numerically more stable. Applying C-RONS to the Kuramoto--Sivashinsky PDE, we observed a 300 times speedup in the computation, while simultaneously reducing the error by two orders of magnitude. Although here we only considered equidistant collocation points, choosing them on an unstructured grid or even an adaptive grid is certainly a possibility.

The RONS equations take the form of a system of nonlinear ODEs which evolve the parameters of the approximate solution. These ODEs tend to become stiff as the number of parameters $n$ increases. As a result, one either has to take exceedingly small time steps or use implicit time integration schemes. To address this issue, we introduced regularized versions of S-RONS and C-RONS that allow fast time integration even with explicit schemes. Regularized RONS adds Tikhonov penalization to the underlying minimization problem such that the resulting ODEs are not stiff. Applying this regularization to both Fokker--Planck and Kuramoto--Sivashinsky equations led to significant speedup without adversely affecting the accuracy of the solutions. Although here we chose the regularization parameter in an ad hoc manner, rigorous methods exist for determining the optimal choice of this parameter~\cite{engl1996,Ito2011}.

The computational methods developed here pave the way for RONS to be used as a shape-morphing spectral method for efficient numerical solution of nonlinear PDEs. In contrast to existing spectral methods where the modes are static in time, the RONS-based spectral methods will allow the modes to change shape and adapt to the solution of the PDE. As a result, these methods will be specially suitable for solving PDEs with localized features (e.g., sharp gradients or shocks) and for advection-dominated PDEs. Future work will explore this avenue by determining the appropriate choice of the shape-morphing modes and carrying out error analysis of the resulting spectral method.

\subsection*{Funding}
This work was supported by the National Science Foundation through the award DMS-2208541.

\appendix

\section{Proof of Theorem~\ref{thm:tikh_soln}}
\label{sec:tikh_proof}

We first note that by taking a time derivative, we can write the constraints in~\eqref{eq:min_cons_tikh} as 
\begin{equation}
	\frac{\id}{\id t}I_k(\vc q(t)) =\langle \nabla I_k(\vc q),\dot{\vc q}\rangle =0, \quad k = 1,2,...,m.
\end{equation}
Introducing the Lagrange multipliers $\hat {\pmb \lambda }= (\hat \lambda_1,...,\hat \lambda_m)^\top \in \R^m$, we define the augmented cost function
\begin{equation}
	\hat{  \mathcal J}_c(\vc q,\dot{\vc q},\hat {\pmb \lambda } ) := \hat{  \mathcal J}(\vc q,\dot{\vc q}) + \sum_{ k = 1 }^{ m } \hat \lambda_k\langle \nabla I_k(\vc q),\dot{\vc q}\rangle.
	\label{eq:AugCostFunc}
\end{equation}	
If a solution to the constrained optimization problem~\eqref{eq:min_cons_tikh} exists, the partial derivatives of $\hat{  \mathcal J}_c$ with respect to $q_i$ and $\hat \lambda_k$ must vanish at the minimizer. This yields
\begin{subequations}
	\begin{equation}
		\nabla_{\dot{\vc q}} \hat{  \mathcal J}+  \sum_{ k = 1 }^{ m } \hat \lambda_k \grad I_k =0,
		\label{eq:augLag_grad1n}
	\end{equation}
	\begin{equation}
		\langle \nabla I_1(\vc q),\dot{\vc q}\rangle = \langle \nabla I_2(\vc q),\dot{\vc q}\rangle= ... = \langle \nabla I_m(\vc q),\dot{\vc q}\rangle = 0.
		\label{eq:augLag_grad2n}
	\end{equation}
\end{subequations}	
We already know (see \cite{anderson2021}, Theorem 1) that our original cost functional satisfies  $\nabla_{\dot{\vc q}}\mathcal J= M(\vc q) \dot{\vc q} - \vc f(\vc q)$. 
Similarly, we can calculate the gradient of  our regularized cost functional to obtain $\nabla_{\dot{\vc q}} \hat{\mathcal J}= (M(\vc q)+\Gamma^\top\Gamma) \dot{\vc q} - \vc f(\vc q)$.

We note that the matrix $(M(\vc q)+\Gamma^\top\Gamma)$ is symmetric positive definite. This is because the metric tensor $M$ is symmetric positive semi-definite, and $\Gamma^\top\Gamma$ is symmetric positive definite due to the assumption that $\Gamma$ is full column rank. Therefore, $(M(\vc q)+\Gamma^\top\Gamma)$ is symmetric positive definite and thus invertible. For notational convenience, we define the regularized metric tensor $\hat M (\vc q) := (M(\vc q)+\Gamma^\top\Gamma)$.

Using the fact that $\hat M$ is invertible,  equation~\eqref{eq:augLag_grad1n} yields
\begin{equation}
	\dot{\vc q} = \hat M ^{-1} (\vc q) \left[ \vc f(\vc q) -\sum_{ k = 1 }^{ m } \hat{\lambda}_k \nabla I_k(\vc q)\right].
	\label{eq:proof_qdot}
\end{equation}
Substituting this expression into~\eqref{eq:augLag_grad2n}, we obtain $m$ equations
\begin{equation}
	\sum_{ k = 1 }^{ m } \hat{\lambda}_k \langle \grad I_i , \hat M^{-1}(\vc q) \grad I_k\rangle = \langle \nabla I_i , \hat M^{-1} (\vc q) \vc f\rangle, \quad i = 1,2,...,m.
	\label{eq:augLag_gradn}
\end{equation}
Equation~\eqref{eq:augLag_gradn} can be written as the system of equations
$\hat C \hat{ \pmb \lambda}  = \hat{ \vc b },$
where $\hat C$ is the \textit{regularized constraint matrix} with entries given by
\begin{equation}
	\hat C_{ij} = \langle \grad I_j ,\hat M^{-1} \grad I_i\rangle,\quad i,j\in\{1,2,\cdots,m\}, 
\end{equation}
and the vector $\hat {\vc b} = (\hat b_1, \hat b_2,\cdots, \hat b_m)^\top\in\R^m$ is given by
\begin{equation}
	\hat b_i = \langle \grad I_i , \hat M^{-1} \vc f \rangle,\quad i=1,2,\cdots, m.
\end{equation}
The matrix $\hat C$ is symmetric positive definite, provided that the constraint gradients $\nabla I_1(\vc q), \nabla I_2(\vc q)$, $\cdots$, $\nabla I_m(\vc q)$ are linearly independent (see \cite{anderson2021}, Lemma 2). 
Thus, the Lagrange multipliers $\hat{ \pmb \lambda} $ are the uniquely determined by $ \hat{ \pmb \lambda}  = \hat C^{-1}\hat{ \vc b }$. 
Therefore, $\dot{\vc q}$ must satisfy equation~\eqref{eq:proof_qdot} where the Lagrange multipliers are determined by solving the system of equations $\hat C \hat{ \pmb \lambda}  = \hat{ \vc b }$.

%\bibliographystyle{plain}
%\bibliography{../bibliog}

\begin{thebibliography}{10}
	
	\bibitem{Adcock12}
	T.~A.~A. Adcock, R.~H. Gibbs, and P.~H. Taylor.
	\newblock The nonlinear evolution and approximate scaling of directionally
	spread wave groups on deep water.
	\newblock {\em Proc. R. Soci. A}, 468(2145):2704--2721, 2012.
	
	\bibitem{adcock09}
	T.~A.~A. Adcock and P.~H. Taylor.
	\newblock Focusing of unidirectional wave groups on deep water: an approximate
	nonlinear {S}chr\"odinger equation-based model.
	\newblock {\em Proc. R. Soci. A}, 465(2110):3083--3102, 2009.
	
	\bibitem{anderson2021}
	W.~Anderson and M.~Farazmand.
	\newblock Evolution of nonlinear reduced-order solutions for {PDE}s with
	conserved quantities.
	\newblock {\em SIAM J. on Scientific Computing}, 44:A176--A197, 2022.
	
	\bibitem{anderson22}
	W.~Anderson and M.~Farazmand.
	\newblock Shape-morphing reduced-order models for nonlinear schrödinger
	equations.
	\newblock {\em Nonlinear Dyn}, 108:2889--2902, 2022.
	
	\bibitem{babaee17}
	H.~Babaee, M.~Farazmand, G.~Haller, and T.~P. Sapsis.
	\newblock Reduced-order description of transient instabilities and computation
	of finite-time {L}yapunov exponents.
	\newblock {\em Chaos}, 27(6):063103, 2017.
	
	\bibitem{otd}
	H.~Babaee and T.~P. Sapsis.
	\newblock A minimization principle for the description of modes associated with
	finite-time instabilities.
	\newblock {\em Proc. R. Soc. A}, 472(2186), 2016.
	
	\bibitem{beale1985}
	J.~T. Beale and A.~Majda.
	\newblock High order accurate vortex methods with explicit velocity kernels.
	\newblock {\em J. Comput. Phys.}, 58(2):188--208, 1985.
	
	\bibitem{Willcox2015}
	P.~Benner, S.~Gugercin, and K.~Willcox.
	\newblock A survey of projection-based model reduction methods for parametric
	dynamical systems.
	\newblock {\em SIAM Review}, 57(4):483--531, 2015.
	
	\bibitem{Bridges2006}
	T.~J. Bridges and S.~Reich.
	\newblock Numerical methods for {H}amiltonian {PDE}s.
	\newblock {\em Journal of Physics A: Mathematical and General}, 39(19):5287,
	apr 2006.
	
	\bibitem{Bruna2023_github}
	J.~Bruna, B.~Peherstorfer, and E.~Vanden-{E}ijnden.
	\newblock Github repository: Neural {G}alerkin with active learning for
	high-dimensional evolution equations.
	\newblock \url{https://github.com/pehersto/ng}, 2022.
	
	\bibitem{bruna22}
	J.~Bruna, B.~Peherstorfer, and E.~Vanden-{E}ijnden.
	\newblock Neural galerkin scheme with active learning for high-dimensional
	evolution equations, 2022.
	
	\bibitem{calvetti2000}
	D.~Calvetti, S.~Morigi, L.~Reichel, and F.~Sgallari.
	\newblock Tikhonov regularization and the l-curve for large discrete ill-posed
	problems.
	\newblock {\em Journal of Computational and Applied Mathematics},
	123(1):423--446, 2000.
	
	\bibitem{carlberg2018}
	K.~Carlberg, Y.~Choi, and S.~Sargsyan.
	\newblock Conservative model reduction for finite-volume models.
	\newblock {\em Journal of Computational Physics}, 371:280--314, 2018.
	
	\bibitem{Cifani2022}
	P.~Cifani, M.~Viviani, E.~Luesink, K.~Modin, and B.~J. Geurts.
	\newblock Casimir preserving spectrum of two-dimensional turbulence.
	\newblock {\em Phys. Rev. Fluids}, 7:L082601, Aug 2022.
	
	\bibitem{koumoutsakos_2000}
	G.-H. Cottet and P.~D. Koumoutsakos.
	\newblock {\em Vortex Methods: Theory and Practice}.
	\newblock Cambridge University Press, 2000.
	
	\bibitem{cousins15}
	W.~Cousins and T.~P. Sapsis.
	\newblock Unsteady evolution of localized unidirectional deep-water wave
	groups.
	\newblock {\em Phys. Rev. E}, 91(6):063204, 2015.
	
	\bibitem{cvitanovic2010}
	P.~Cvitanovi\'{c}, R.~L. Davidchack, and E.~Siminos.
	\newblock On the state space geometry of the {K}uramoto--{S}ivashinsky flow in
	a periodic domain.
	\newblock {\em SIAM Journal on Applied Dynamical Systems}, 9(1):1--33, 2010.
	
	\bibitem{Babaee2022}
	M.~Donello, M.~H. Carpenter, and H.~Babaee.
	\newblock Computing sensitivities in evolutionary systems: {A} real-time
	reduced order modeling strategy.
	\newblock {\em SIAM Journal on Scientific Computing}, 44(1):A128--A149, 2022.
	
	\bibitem{dong21}
	S.~Dong and N.~Ni.
	\newblock A method for representing periodic functions and enforcing exactly
	periodic boundary conditions with deep neural networks.
	\newblock {\em Journal of Computational Physics}, 435:110242, 2021.
	
	\bibitem{dormand1980}
	J.R. Dormand and P.J. Prince.
	\newblock A family of embedded runge-kutta formulae.
	\newblock {\em Journal of Computational and Applied Mathematics}, 6(1):19--26,
	1980.
	
	\bibitem{du21}
	Y.~Du and T.A. Zaki.
	\newblock Evolutional deep neural network.
	\newblock {\em Phys. Rev. E}, 104:045303, Oct 2021.
	
	\bibitem{engl1996}
	H.~W. Engl, M.~Hanke, and A.~Neubauer.
	\newblock {\em Regularization of inverse problems}, volume 375.
	\newblock Springer Science \& Business Media, The Netherlands, 1996.
	
	\bibitem{faraz_adjoint}
	M.~Farazmand.
	\newblock An adjoint-based approach for finding invariant solutions of
	{N}avier-{S}tokes equations.
	\newblock {\em J. Fluid Mech.}, 795:278--312, 2016.
	
	\bibitem{PRE2016}
	M.~Farazmand and T.~P. Sapsis.
	\newblock Dynamical indicators for the prediction of bursting phenomena in
	high-dimensional systems.
	\newblock {\em Phys. Rev. E}, 94:032212, 2016.
	
	\bibitem{Gill2021}
	P.~E. Gill, W.~Murray, and M.~H. Wright.
	\newblock {\em Numerical Linear Algebra and Optimization}.
	\newblock Society for Industrial and Applied Mathematics, Philadelphia, PA,
	2021.
	
	\bibitem{golub1999}
	G.H. Golub, P.C. Hansen, and D.P. O'Leary.
	\newblock Tikhonov regularization and total least squares.
	\newblock {\em SIAM Journal on Matrix Analysis and Applications},
	21(1):185--194, 1999.
	
	\bibitem{Weinan2018}
	J.~Han, A.~Jentzen, and W.~E.
	\newblock Solving high-dimensional partial differential equations using deep
	learning.
	\newblock {\em Proceedings of the National Academy of Sciences},
	115(34):8505--8510, 2018.
	
	\bibitem{hyman1986}
	J.~M. Hyman and B.~Nicolaenko.
	\newblock The {K}uramoto--{S}ivashinsky equation: {A} bridge between {PDE}'s
	and dynamical systems.
	\newblock {\em Physica D: Nonlinear Phenomena}, 18(1):113--126, 1986.
	
	\bibitem{Ipsen2009}
	I.~C.~F. Ipsen.
	\newblock {\em Numerical Matrix Analysis}.
	\newblock Society for Industrial and Applied Mathematics, Philadelphia, USA,
	2009.
	
	\bibitem{Ito2011}
	K.~Ito, B.~Jin, and T.~Takeuchi.
	\newblock A regularization parameter for nonsmooth {T}ikhonov regularization.
	\newblock {\em SIAM Journal on Scientific Computing}, 33(3):1415--1438, 2011.
	
	\bibitem{karniadakis2005}
	G.~E. Karniadakis and S.~Sherwin.
	\newblock {\em Spectral/hp element methods for computational fluid dynamics}.
	\newblock Oxford University Press, Oxford, UK, 2005.
	
	\bibitem{Keller1992}
	P.~S. Keller.
	\newblock Chaotic behavior of {N}ewton's method.
	\newblock {\em Real Analysis Exchange}, 18(2):490--507, 1992.
	
	\bibitem{Kevrekidis1990}
	I.~G. Kevrekidis, B.~Nicolaenko, and J.~C. Scovel.
	\newblock Back in the saddle again: {A} computer assisted study of the
	{K}uramoto--{S}ivashinsky equation.
	\newblock {\em SIAM Journal on Applied Mathematics}, 50(3):760--790, 1990.
	
	\bibitem{majda2012}
	A.~J. Majda and Y.~Yuan.
	\newblock Fundamental limitations of ad hoc linear and quadratic multi-level
	regression models for physical systems.
	\newblock {\em Discrete \& Continuous Dynamical Systems - B},
	17(1531-3492\_2012\_4\_1333):1333, 2012.
	
	\bibitem{McLachlan2006}
	R.~I. McLachlan and G.~R.~W. Quispel.
	\newblock Geometric integrators for {ODE}s.
	\newblock {\em Journal of Physics A: Mathematical and General}, 39(19):5251,
	apr 2006.
	
	\bibitem{newton_Nvortex}
	P.~K. Newton.
	\newblock {\em The N-vortex problem: analytical techniques}, volume 145 of {\em
		Applied Mathematical Sciences}.
	\newblock Springer, 2001.
	
	\bibitem{pathak2018}
	J.~Pathak, B.~Hunt, M.~Girvan, Z.~Lu, and E.~Ott.
	\newblock Model-free prediction of large spatiotemporally chaotic systems from
	data: A reservoir computing approach.
	\newblock {\em Phys. Rev. Lett.}, 120:024102, 2018.
	
	\bibitem{peng16}
	L.~Peng and K.~Mohseni.
	\newblock Symplectic model reduction of hamiltonian systems.
	\newblock {\em SIAM Journal on Scientific Computing}, 38(1):A1--A27, 2016.
	
	\bibitem{PerezGarcia1996}
	V.~M. P\'erez-Garc\'{\i}a, H.~Michinel, J.~I. Cirac, M.~Lewenstein, and
	P.~Zoller.
	\newblock Low energy excitations of a {B}ose-{E}instein condensate: {A}
	time-dependent variational analysis.
	\newblock {\em Phys. Rev. Lett.}, 77:5320--5323, Dec 1996.
	
	\bibitem{rowley2017}
	C.~W. Rowley and S.~T.~M. Dawson.
	\newblock Model reduction for flow analysis and control.
	\newblock {\em Annual Review of Fluid Mechanics}, 49(1):387--417, 2017.
	
	\bibitem{ruban2015}
	V.~P. Ruban.
	\newblock Anomalous wave as a result of the collision of two wave groups on the
	sea surface.
	\newblock {\em JETP Letters}, 102(10):650--654, 2015.
	
	\bibitem{ruban2015b}
	V.~P. Ruban.
	\newblock Gaussian variational ansatz in the problem of anomalous sea waves:
	{C}omparison with direct numerical simulation.
	\newblock {\em Journal of Experimental and Theoretical Physics},
	120(5):925--932, 2015.
	
	\bibitem{shampine1997}
	L.~F. Shampine and M.~W. Reichelt.
	\newblock The matlab ode suite.
	\newblock {\em SIAM Journal on Scientific Computing}, 18(1):1--22, 1997.
	
	\bibitem{Spiliopoulos2018}
	J.~Sirignano and K.~Spiliopoulos.
	\newblock {DGM}: {A} deep learning algorithm for solving partial differential
	equations.
	\newblock {\em Journal of Computational Physics}, 375:1339--1364, 2018.
	
	\bibitem{yazdani2020}
	A.~Yazdani, L.~Lu, M.~Raissi, and G.~E. Karniadakis.
	\newblock Systems biology informed deep learning for inferring parameters and
	hidden dynamics.
	\newblock {\em PLOS Computational Biology}, 16(11):1--19, 11 2020.
	
\end{thebibliography}

\end{document}